\documentstyle{amsppt}
\pagewidth{15 cm}
\pageheight{22 cm}
\voffset 1 cm
\hoffset 0.45 cm
\TagsOnRight
\topmatter
\title $Tb$ theorem on product spaces \endtitle
\leftheadtext{}
\author Yongsheng Han, Ming-Yi Lee, and Chin-Cheng Lin \endauthor
\rightheadtext{Y. Han, M.-Y. Lee, and C.-C. Lin}
\leftheadtext{$Tb$ theorem on product spaces}
\address  Department of Mathematics,
          Auburn University,
          Auburn, Alabama 36849-5310, U.S.A.  \endaddress
\email hanyong\@mail.auburn.edu \endemail
\address  Department of Mathematics,
          National Central University,
          Chung-Li, Taiwan 320,
          Republic of \break China \endaddress
\email mylee\@math.ncu.edu.tw;\ clin\@math.ncu.edu.tw \endemail
\keywords Almost orthogonality, Carleson measure, Journ\'e's class, Littlewood-Paley estimate, para-accretive function, singular integral operator \endkeywords
\subjclass 42B30 \endsubjclass
\thanks Research by the second and third authors supported by NSC of Taiwan under Grant
        \#NSC 99-2115-M-008-002-MY3 and Grant \#NSC 100-2115-M-008-002-MY3, respectively. \endthanks
\abstract
In this paper, we prove a $Tb$ theorem on product spaces $\Bbb R^n\times \Bbb R^m$,
where $b(x_1,x_2)=b_1(x_1)b_2(x_2)$, $b_1$ and $b_2$ are
para-accretive functions on $\Bbb R^n$ and $\Bbb R^m$, respectively.
\endabstract
\endtopmatter
\document

\head \S1. Introduction
\endhead
In their well-known theory of singular integral operators,
Caldern\'on and Zygmund established the $L^p, 1<p<\infty,$
boundedness of certain convolution singular integral operators on
$\Bbb R^n,$ which generalize the Hilbert transform on $\Bbb R^1$.
This theory has been generalized in two ways: First, the convolution
singular integral operators were replaced by non-convolution
singular integral operators. To be more precise, a continuous
complex-valued function $k(x,y)$ defined on  $\Bbb R^n\times \Bbb
R^n\backslash \{ (x, y): x=y\}$ is called a {\it Calder\'on-Zygmund
kernel} if there exist constant $C>0$ and a regularity exponent
$\varepsilon\in (0,1]$ such that \roster
\item"(i)" $|k(x,y)|\le C|x-y|^{-n}$ \item"(ii)"
$|k(x,y)-k(x',y)|\le C |x-x'|^\varepsilon|x-y|^{-n-\varepsilon}$
             \qquad if\ \ $|x-x'|\le |x-y|/2$\hfill(1.1)
\item"(iii)" $|k(x,y)-k(x,y')|\le C |y-y'|^\varepsilon|x-y|^{-n-\varepsilon}$
             \qquad if\ \ $|y-y'|\le |x-y|/2$.
\endroster
The smallest such constant $C$ is denoted by $|k|_{CZ}.$ We say that
an operator $T$ is a classical singular integral operator if the
operator $T$ is a continuous linear operator from $C^\infty_0(\Bbb
R^n)$ into its dual associated with a Calder\'on-Zygmund kernel
$k(x,y)$ given by
$$\langle Tf, g \rangle=\iint g(x)k(x,y)f(y)dydx$$
for all functions $f, g\in C^\infty_0(\Bbb R^n)$ with disjoint
supports. $T$ is said to be a Calder\'on-Zygmund operator
if it extends to be a bounded operator on $L^2(\Bbb R^n).$
If $T$ is a Calder\'on-Zygmund operator associated with a kernel $k$,
its operator norm is defined by
$\|T\|_{CZ}=\|T\|_{L^2\mapsto L^2}+ |k|_{CZ}$.
 Of course, in general, one cannot conclude that
  such a Calder\'on-Zygmund singular integral operator
 $T$ is bounded on $L^2(\Bbb R^n)$ because
 Plancherel's theorem doesn't work for non-convolution operators.
 However, if one assumes that $T$ is bounded on $L^2(\Bbb R^n),$
 then the $L^p, 1<p<\infty,$ boundedness follows from
 Caldern\'on-Zygmund's real variable method. The $L^2(\Bbb R^n)$
 boundedness of non-convolution singular integral operators was finally
 proved by the remarkable $T1$ theorem of David and Journ\'e in
 [DJ],which gives a general criterion for the $L^2$-boundedness
  of Calder\'on-Zygmund singular integral operators. Unfortunately,
  the $T1$ theorem cannot be applied to
  the Cauchy integral on a Lipschitz curve defined by
$$C(f)(x)=\frac 1\pi\ \text{p.v.}\int_{-\infty}^\infty \frac{f(y)}{(x-y)+i(a(x)-a(y))} dy,$$
where the function $a(x)$ satisfies the Lipschitz condition.

Indeed, it is still an open problem that, without assuming the $L^2$-boundedness,
one does not know how to prove that the Cauchy integral $C(f)$
  on a Lipschitz curve maps the function $1$ into a $BMO$ function.
  Meyer first observed that $C(b)=0$ provided $b(x)=1+ia'(x)$. Therefore, if the function $1$ in
  the $T1$ theorem is allowed to be replaced by an accretive function $b$ which is a bounded complex-valued
  function satisfying $\text{Re}\, b(x)\ge \delta >0$ almost everywhere,  then this result would imply the
  $L^2$-boundedness of the Cauchy integrals on all Lipschitz curves.
  McIntosh and Meyer [MM] obtained such a $Tb$ theorem; that is,
  the $T1$ theorem still holds while the function $1$ is replaced by an accretive function $b$.
  Finally, David, Journ\'e and Semmes [DJS] proved a new $Tb$ theorem
  which says that the function $1$ in the $T1$
  theorem can be replaced by the so-called para-accretive functions $b$ (see definition below).
  See [MM] and [DJS] for more details about the $Tb$ theorems.

 Secondly, by taking the space $\Bbb R^n\times \Bbb R^m$ along with
 two parameter family of dilations $(x,y)\mapsto (\delta_1
 x,\delta_2 y), x\in \Bbb R^n, y\in \Bbb R^m, \delta_i>0, i=1, 2,$
 instead of the classical one-parameter dilation, R. Fefferman and
 Stein [FS] studied the product convolution singular integral operators
 which satisfy analogous conditions enjoyed by the double Hilbert
 transform defined on $\Bbb R\times\Bbb R.$ Journ\'e [J] generalized
 the product convolution singular integral operators to the product
 non-convolution singular integral  operators and introduced a class of
 singular integral operators
 which coincides with the product convolution
singular integral operators with two parameters. More precisely, a
singular integral operator $T$ is said to be in {\it Journ\'e's class}  if
$T$ is a continuous linear operator from $C^\infty_0(\Bbb R)\otimes
C^\infty_0(\Bbb R)\rightarrow[C^\infty_0(\Bbb R)\otimes
C^\infty_0(\Bbb R)]'$ defined by
$$\langle g\otimes k, Tf\otimes h\rangle
   =\iint g(x)\langle k, K_1(x,y)h\rangle f(y) dxdy$$
and
$$\langle k\otimes g, Th\otimes f\rangle =\iint g(x)\langle k, K_2(x,y)h\rangle f(y) dxdy,$$
for all $f, g, h, k\in C^\infty_0(\Bbb R)$ with supp$(f)\ \cap$ supp$(g)=\varnothing$,
and a pair $(K_1,K_2)$ of $\delta CZ$-$\delta$-standard kernels
defined in [J, p.63]. Moreover, Journ\'e [J] proved the product $T1$ theorem as follows.

\proclaim{Theorem A} Let $T$ belong to Journ\'e's class. Then $T$ and $\widetilde T$ are
bounded on $L^2(\Bbb R^2)$ if and only if $T1$, ${}^tT1$,
$\widetilde T1$, and ${}^t\widetilde T1$ lie in $BMO(\Bbb R\times
\Bbb R)$ and $T$ has the weak boundedness property.
\endproclaim
Here, ${}^tT$ is the transport of $T$ and $\widetilde T,$ the
partial adjoint operator of $T,$ is defined by $\langle g\otimes k,
\widetilde Tf\otimes h\rangle =\langle f\otimes k, Tg\otimes h
\rangle.$

The purpose of this paper is to unify up to a certain
generalizations of the $Tb$ theorem in [DJS], the product $T1$
theorem in [J], and the product $Tb$ theorem with $Tb={}^tTb=\widetilde Tb={}^t\widetilde Tb=0$ in [LZ].
In order to state our main result, the $Tb$ theorem
on the product space, we first recall some basic definitions and
notations.

Let $C^\eta_0(\Bbb R^n)$ denote the space of continuous functions
$f$ with compact support such that
$$\sup\limits_{x\ne y} \frac{|f(x)-f(y)|}{|x-y|^\eta}<\infty$$
and let $C^\eta_0(\Bbb R^n\times \Bbb R^m), \eta>0$,
  denote the space of continuous functions $f$ with compact support such that
$$\| f\|_{\eta} : = \sup_{\Sb x_1\ne y_1\\ x_2\ne y_2\endSb}
  \frac{|f(x_1,x_2)-f(y_1,x_2)-f(x_1,y_2)+f(y_1,y_2)|}{|x_1-y_1|^\eta |x_2-y_2|^\eta} <\infty.$$
A {\it singular integral operator} $T$ is a continuous linear
operator from $C^\eta_0(\Bbb R^n\times \Bbb R^m)$  into its dual
$(C^\eta_0(\Bbb R^n\times \Bbb R^m))'$ associated with a kernel $K(x_1,x_2,y_1,y_2)$,
 a continuous complex-valued function on
$\Bbb R^n\times \Bbb R^m \times \Bbb R^n\times \Bbb R^m\backslash \{(x_1,x_2,y_1,y_2): x_1=y_1\ \text{or}\ x_2=y_2\}$,
 and it can be defined on $C^\eta_0(\Bbb R^n)\otimes C^\eta_0(\Bbb R^m)$ as follows
$$\langle Tf_1\otimes f_2, g_1\otimes g_2\rangle
= \int_{\Bbb R^n\times \Bbb R^m}\int_{\Bbb R^n\times \Bbb R^m}
    g_1(x_1)g_2(x_2)K(x_1,x_2,y_1,y_2)f_1(y_1)f_2(y_2)dx_1dx_2dy_1dy_2    $$
    for all $f_1,g_1\in C^\eta_0(\Bbb R^n)$ with supp$f_1\
     \cap$ supp$g_1=\varnothing$ and $f_2,g_2\in C^\eta_0(\Bbb R^m)$
with supp$f_2\ \cap$ supp$g_2=\varnothing,$ where
$K(x_1,x_2,y_1,y_2),$ the kernel of $T,$ satisfies the following
conditions: for each $x_1,y_1\in \Bbb R^n$,
  $\widetilde K^1(x_1,y_1)$ is a Calder\'on-Zygmund operator acting on functions
   on $\Bbb R^m$ with the kernel \break $\widetilde K^1(x_1,y_1)(x_2,y_2)=K(x_1,x_2,y_1,y_2),$
   and similarly, for each
   $x_2,y_2\in \Bbb R^m$, $\widetilde K^2(x_2,y_2)$ is a
Calder\'on-Zygmund operator acting on functions on $\Bbb R^n$ with
the kernel $\widetilde K^2(x_2,y_2)(x_1,y_1)=K(x_1,x_2,y_1,y_2).$
Moreover, there exist constants $C>0$ and
$\varepsilon\in (0,1]$ such that \roster
\item"($A_1$)" $\big\|\widetilde K^1(x_1,y_1)\big\|_{CZ}\le C|x_1-y_1|^{-n}$,
\vskip 0.1cm
\item "" $\big\|\widetilde K^1(x_1,y_1)-\widetilde K^1(x_1,y'_1)\big\|_{CZ}
                \le C |y_1-y'_1|^\varepsilon|x_1-y_1|^{-(n+\varepsilon)}\qquad\text{for}\ \
                    |y_1-y'_1|\le |x_1-y_1|/2,$
\vskip 0.1cm
\item""  $\big\|\widetilde K^1(x_1,y_1)-\widetilde K^1(x'_1,y_1)\big\|_{CZ}
                \le C |x_1-x'_1|^\varepsilon|x_1-y_1|^{-(n+\varepsilon)}\qquad\text{for}\ \
                    |x_1-x'_1|\le |x_1-y_1|/2,$
\vskip 0.1cm
\item"($A_2$)"$\big\|\widetilde K^2(x_2,y_2)\big\|_{CZ}\le C|x_2-y_2|^{-m}$ ,
\vskip 0.1cm
\item "" $\big\|\widetilde K^2(x_2,y_2)-\widetilde K^2(x_2,y'_2)\big\|_{CZ}
           \le C |y_2-y'_2|^\varepsilon|x_2-y_2|^{-(m+\varepsilon)}\qquad\text{for}\ \
                     |y_2-y_2'|\le |x_2-y_2|/2,$
\vskip 0.1cm
\item "" $\big\|\widetilde K^2(x_2,y_2)-\tilde K^2(x'_2,y_2)\big\|_{CZ}
           \le C
           |x_2-x'_2|^\varepsilon|x_2-y_2|^{-(m+\varepsilon)}\qquad\text{for}\
           \     |x_2-x_2'|\le |x_2-y_2|/2.$
\endroster
Let $T$ be a singular integral operator. For $f_1, g_1\in C^\eta_0 (\Bbb R^n)$ and $f_2,g_2\in C^\eta_0(\Bbb R^m)$,
the operator $\langle T^1f_1, g_1 \rangle : C^\eta_0 (\Bbb R^m) \mapsto (C^\eta_0 (\Bbb R^m))'$ is defined by
$$\big\langle \langle T^1f_1, g_1 \rangle f_2, g_2\big\rangle= \langle Tf_1\otimes f_2, g_1\otimes g_2\rangle.$$
It is easy to see that $\langle T^1f_1, g_1 \rangle$ is a singular integral operator on $\Bbb R^m$
with kernel $\langle T^1f_1, g_1 \rangle(x_2,y_2)= \langle \widetilde K^2(x_2,y_2)f_1, g_1\rangle$.
One defines $\langle T^2f_2, g_2 \rangle : C^\eta_0 (\Bbb R^n) \mapsto (C^\eta_0 (\Bbb R^n))'$ similarly.
We say that $T$ has the {\it weak boundedness property}, denoted by $T\in WBP$, if there exists
$C>0$ such that for all functions $f_1, g_1\in C^\eta_0 (\Bbb R^n)$
supported in a cube $Q_1$ and $f_2,g_2\in C^\eta_0(\Bbb R^m)$
supported in another cube $Q_2$,
$$\cases  \|\langle T^1f_1, g_1 \rangle\|_{CZ}\le C|Q_1|^{1+\frac{2\eta} n}
\|f_1\|_{\eta(\Bbb R^n)}\|g_1\|_{\eta(\Bbb R^n)}    \\
      \|\langle T^2f_2, g_2 \rangle\|_{CZ}\le C|Q_2|^{1+\frac {2\eta} m}
\|f_2\|_{\eta(\Bbb R^m)}\|g_2\|_{\eta(\Bbb R^m)} \endcases.$$

A bounded complex-valued function $b$ defined on $\Bbb R^n$ is said
to be {\it para-accretive} if there exist constants $C, \gamma >0$
such that, for all cubes $Q\subset\Bbb R^n$, there is a
$Q^{\prime}\subset Q$ with $\gamma |Q|\le |Q^{\prime}|$ satisfying
$$\frac{1}{|Q|}\Big|\int_{Q^{\prime}} b(x)\,dx\Big|\ge C.$$
Note that, by the Lebesgue differentiation theorem, $b^{-1}(x)$ is
also bounded.

Suppose $b(x_1,x_2)=b_1(x_1)b_2(x_2),$ where $b_1$ and $b_2$ are para-accretive functions on
$\Bbb R^n$ and $\Bbb R^m,$ respectively. A {\it
generalized singular integral operator} is a continuous linear
operator $T$ from $bC^\eta_0(\Bbb R^n \times \Bbb R^m)$ into $(bC^\eta_0(\Bbb R^n \times \Bbb R^m))'$
for all $\eta>0$ if the kernel of $T$ is a singular integral kernel and
 for $f_1,f_2,g_1,g_2\in C^\eta_0(\Bbb R^n)$ with
supp$(f_1)\ \cap$ supp$(g_1)=\varnothing$ and $f_2,g_2\in C^\eta_0(\Bbb R^m)$
with supp$(f_2)\ \cap$ supp$(g_2)=\varnothing,$
$$\align
&\langle M_bTM_bf_1\otimes f_2, g_1\otimes g_2\rangle
= \int_{\Bbb R^n\times \Bbb R^m}\int_{\Bbb R^n\times \Bbb R^m}
b_2(x_2)b_1(x_1)g_1(x_1)g_2(x_2) \\
&\hskip 6cm\times K(x_1,x_2,y_1,y_2)b_2(y_2)b_1(y_1)f_1(y_1)f_2(y_2)dx_1dx_2dy_1dy_2,
\endalign$$
where $M_b$ denotes the multiplication operator by $b$; that is, $M_bf=bf$.

Suppose that $T$ is a generalized singular integral operator associated to a kernel
$K(x_1,x_2,y_1,y_2)$. Then ${}^tT$, $\widetilde T$, and
${}^t\widetilde T$ are singular integral operators associated to
kernels ${}^t\!K(x_1,x_2,y_1,y_2):=K(y_1,y_2,x_1,x_2)$,
$\widetilde K(x_1,x_2,y_1,y_2):=K(y_1,x_2,x_1,y_2)$,
and ${}^t\!\widetilde K(x_1,x_2,y_1,y_2):=K(x_1,y_2,y_1,x_2)$, respectively.
Our main result is the following

\proclaim{Theorem 1} Suppose that $b_1$ and $b_2$ are para-accretive
functions on $\Bbb R^n$ and $\Bbb R^m$, respectively,
$b(x_1,x_2)=b_1(x_1)b_2(x_2)$, and $T$ is a generalized singular
integral operator. If $Tb$, ${}^tTb$, $\widetilde Tb$,
${}^t\widetilde Tb \in BMO(\Bbb R^n\times \Bbb R^m)$ and $M_bTM_b\in WBP$,
then $T$ is bounded on $L^2(\Bbb R^{n+m})$
\endproclaim

Applying the above Theorem 1 together with [J, Theorem 3],
we obtain the following product $Tb$ theorem.

\proclaim{Theorem 2} Suppose that $b_1$ and $b_2$ are para-accretive
functions on $\Bbb R^n$ and $\Bbb R^m$, respectively,
$b(x_1,x_2)=b_1(x_1)b_2(x_2)$, and $T$ is a generalized singular
integral operator. Then $T$ and $\widetilde T$ are bounded on
$L^2(\Bbb R^{n+m})$ if and only if $Tb$, ${}^tTb$, $\widetilde Tb$,
${}^t\widetilde Tb \in BMO(\Bbb R^n\times \Bbb R^m)$ and $M_bTM_b\in WBP$.
\endproclaim
In order to describe our approach to the proof, we first recall the
general philosophy of the proofs of the $T1$ theorem of David and
Journ\'e and the product $T1$ theorem of Journ\'e. The $T1$ theorem
is proved by two steps. In the first step, one considers the case
where $T(1)=T^*(1)=0$ and then uses the Cotlar-Stein lemma.
To be more precise, let $\psi\in C^\infty_0(\Bbb R^n)$ with
$\int^\infty_0 |\widehat \psi(t\xi)|^2\frac{dt}{t}=1$ for all
$\xi\not=0.$ Let $\psi_t(x)=t^{-n}\psi(\frac{x}{t}).$ The operator
$U_j=\int^j_{1/j}\psi_t*\psi_t\frac{dt}{t}$ converges
strongly to the identity on $L^2(\Bbb R^n)$ as $j\rightarrow
\infty.$ Since $U_j$ is continuous on $C^\infty_0(\Bbb R^n), U_jTU_{j'}$ is well defined on
$C^\infty_0(\Bbb R^n)$ for all $j$ and $j'.$ Therefore,
$T$ is bounded on $L^2(\Bbb R^n)$ if and only if $U_jTU_{j'}$ is
bounded on $L^2(\Bbb R^n)$ with a norm independent of $j$ and $j'.$
The second step is to
use the para-product operator to reduce the general case to the
first step. The para-product operator is defined by
$\Pi_b(f)=\int^\infty_0\psi_t*(\psi_t*b(\cdot)\phi_t*f(\cdot))(x)\frac{dt}{t},$
where $b\in BMO$ and $\phi\in C^\infty_0(\Bbb R^n)$ with
$\int\phi(x)dx=1.$ Using a result of the Carleson measure, it was well-known that
the para-product operator $\Pi_b$ is a Calder\'on-Zygmund
operator on $\Bbb R^n,$ moreover, $\Pi_b(1)=b$ and $\Pi^*_b(1)=0.$
We would like to remark that the fact that for any $b\in BMO(\Bbb R^n), \Pi_b(1)=b$ in the sense that
$$\langle f, b\rangle =\Big\langle f, \int^\infty_0\psi_t*\psi_t*b\frac{dt}{t}\Big\rangle $$
for all $f\in H^1(\Bbb R^n).$ Now one can decompose $T$ by $ T={\widetilde T}+
\Pi_{T1}+\Pi^*_{T^*1},$ where ${\widetilde T}=T-\Pi_{T1}-\Pi^*_{T^*1}.$
By the first step and properties of the para-product operators, ${\widetilde T}$ is
bounded on $L^2(\Bbb R^n)$ and hence  the $L^2$ boundedness of $T$
follows. The proof of the product $T1$ theorem follows from a similar
way. In the first step, the $L^2$ boundedness follows from
the assumptions that $T(1)={}^tT(1)=\widetilde T(1)={}^t\widetilde T(1)=0.$
The product-type paraproduct operators are constructed in the
second step and the general case is then reduced to the first step.
See [DJ] and [J] for the details.

In this paper, we will employ a new approach to prove the product
$Tb$ theorem. The new feature of our approach is to use the almost
orthogonality argument to obtain a new decomposition of $T.$ The
para-product operators constructed in [DJ] and [J] are avoided. To
see how this approach works, we outline first a new proof of the
classical $T1$ theorem based on the almost orthogonality argument
and our new decomposition. Here the almost orthogonality argument means
that for the function $\psi$ as given above, there exists a constant $C$ such
that
$$ |\psi_t\ast \psi_s(x)|\leq
C\Big(\frac{t}{s}\wedge\frac{s}{t}\Big)\frac{(t\vee s)}{((t\vee s)
+|x|)^{(n+1)}},$$ where $a \wedge b = \min\{a, b\}$ and $a\vee b=\max
\{a,b\}.$

If $T$ satisfies the cancellation conditions $T(1)=T^*(1)=0,$ then
one still has the following almost orthogonality argument:
$$ \bigg|\int_{\Bbb R^n}\int_{\Bbb R^n}\psi_t(x-u) \ k(u,v) \ \psi_s(v-y) \ du dv\bigg|\leq
    C\Big(\frac{t}{s}\wedge\frac{s}{t}\Big)^{\varepsilon'}
    \frac{(t\vee s)^{\varepsilon'}}{((t\vee s) +|x-y|)^{(n+\varepsilon')}}  ,$$ where
$0<\varepsilon'<\varepsilon$ and $\varepsilon$ is the regularity exponent of
the kernel $k$ given in (1.1).

In general, the above almost orthogonality argument doesn't hold
without the assumptions on the cancellation conditions on the kernel
of $T.$ However, if $t\leq s,$ then
$$\bigg|\int_{\Bbb R^n}\int_{\Bbb R^n}\psi_t(x-u)k(u,v)\big[\psi_s(v-y)-\psi_s(x-y)\big]dudv\bigg|\leq
   C\Big(\frac{t}{s}\Big)^{\varepsilon'}\frac{s^{\varepsilon'}}{(s +|x-y|)^{(n+\varepsilon')}}.$$
This leads to the following decomposition of $T.$ Now
suppose that $T$ satisfies the hypotheses of the $T1$ theorem of
David and Journ\'e. As in the proof of the $T1$ theorem, we would
like to show that $U_jTU_{j'}$ is bounded on $L^2(\Bbb R^n)$ as both $j, j'\to \infty$.
To do this, we decompose the kernel of $\lim_{j,j'\to \infty} U_jTU_{j'}$ as follows,
$$\align
&\lim_{\Sb j\to \infty \\ j'\to \infty \endSb} U_jTU_{j'}(x,y) \\
&=\iiint_{0<t<\infty}\int_{0<s<\infty}\psi_t*\psi_t(x-u)k(u,v)
   \psi_s*\psi_s(v-y)dudv\frac{dt}{t}\frac{ds}{s}\\
&=\iint\hskip-0,2cm\iint\hskip-0,2cm\int_{\Sb s<t<\infty\\
0<s<\infty\endSb}\psi_t(x-u')\big[\psi_t(u'-u)-\psi_t(u'-v')\big]k(u,v)
   \psi_s(v-v')\psi_s(v'-y)du'dudvdv'\frac{dt}{t}\frac{ds}{s}\\
&\quad +\iiint\hskip-0,2cm\iint_{\Sb 0<t<\infty\\ t\le s<\infty\endSb}
    \hskip-0,2cm\psi_t(x-u')\psi_t(u'-u)k(u,v)
   \big[\psi_s(v-v')-\psi_s(u'-v')\big]\psi_s(v'-y)du'dudvdv'\frac{dt}{t}\frac{ds}{s}\\
&\quad +\iint_{\Sb s<t<\infty\\ 0<s<\infty\endSb}\psi_t*\psi_t(x-v')
   \psi_s*T^*(1)(v')\psi_s(v'-y)dv'\frac{dt}{t}\frac{ds}{s}\\
&\quad
+\iint_{\Sb 0<t<\infty\\ t\le s<\infty\endSb}\psi_t(x-u')\psi_t*T(1)(u')
    \psi_s*\psi_s(u'-y)du'\frac{dt}{t}\frac{ds}{s}\\
&:=k_1(x,y)+k_2(x,y)+k_3(x,y)+k_4(x,y).
\endalign$$
Let $f, g\in C^\eta_0(\Bbb R^n).$ The almost orthogonality
argument yields $|\langle g,\langle k_1, f\rangle \rangle |\leq C\|f\|_{L^2}\|g\|_{L^2}$ and similarly
for $k_2.$ The $L^2$ boundedness of operators with the kernels $k_3$
and $k_4$ follows from a result of the Carleson measure.

To carry out the above approach to the proof of the product $Tb$
theorem, we need some definitions and notations.

Let $b$ be a para-accretive function defined on $\Bbb R^n$.
   A sequence of operators $\{S_j\}_{j\in \Bbb Z}$ is called to
   be an {\it approximation to the identity associated to $b$}
   if $S_j(x,y),$ the kernels of $S_j,$ are functions
   from $\Bbb R^n\times \Bbb R^n$
   into $\Bbb C$ such that there exist  constant $C$
   and some $0<\varepsilon\le 1,$ and
   for all $j\in \Bbb Z, x, x^\prime, y$,
   $y^\prime \in {\Bbb R^n}$,
\roster
\item"(i)" $|S_j(x,y)|\le C2^{jn}$ \qquad for all $x,y \in \Bbb R^n$ \vskip0.1cm
\item"(ii)"
$|S_j(x, y)|=0$\qquad if $|x-y|\ge C2^{-j}$,\vskip0.1cm
\item"(iii)"
  $|S_j(x,y) - S_j(x',y)| \le C2^{j(n+\varepsilon)}|x- x^\prime|^\varepsilon$,
\vskip0.1cm
\item"(iv)"
  $|S_j(x, y) - S_j(x, y^\prime) |
    \le C2^{j(n+\varepsilon)}|y- y^\prime|^\varepsilon$,
\vskip0.1cm
\item"(v)"
   $\big|\big[S_j(x, y) - S_j(x, y^\prime)\big] - \big[S_j(x^\prime, y)
   - S_j(x^\prime , y^\prime)\big] \big|$
     $\le C2^{j(n+2\varepsilon)}|x- x^\prime|^\varepsilon|y- y^\prime|^\varepsilon$,
\vskip0.1cm
\item"(vi)"
  $\displaystyle\int_{\Bbb R^n} S_j(x, y) b(y)dy = 1$
  \qquad for all $j\in\Bbb Z$ and $x\in\Bbb R^n$,
\vskip0.1cm
\item"(vii)"
  $\displaystyle\int_{\Bbb R^n} S_j(x, y) b(x)dx = 1$
  \qquad for all $j\in\Bbb Z$ and $y\in\Bbb R^n$.
\endroster

We remark that the existence of such an approximation to the
identity follows from Coifman's idea, which was constructed in
[DJS]. Suppose that $b_1, b_2$ are para-accretive functions on $\Bbb R^n, \Bbb R^m$ respectively,
and $\{S_j\}$, $\{\overline S_k\}$ are
approximations to the identity associated to $b_1, b_2$
respectively. Set $D_j=S_j-S_{j-1}$ and $\overline D_k=\overline S_k-\overline S_{k-1}$.

We recall the {\it $\Cal G$-function associated to a para-accretive function $b$} given by
$$\Cal G(f)(x_1,x_2):=\bigg\{\sum_{j,k\in \Bbb Z} |D_j\overline D_kM_bf(x_1,x_2)|^2\bigg\}^{1/2}
  \qquad\text{for $f\in L^2(\Bbb R^n\times \Bbb R^m)$.}$$
It is known that $\big\|\sum_{j\in \Bbb Z}
D_jM_{b_1}f\big\|_{L^2(\Bbb R^n)}\le C\|f\|_{L^2(\Bbb R^n)}$
and $\big\|\sum_{j\in \Bbb Z}D^N_jM_{b_1}f\big\|_{L^2(\Bbb R^n)}\le C\|f\|_{L^2(\Bbb R^n)}$
for single parameter (see [DJS, p. 19]).
Thus $\Cal G$-function is bounded on $L^2(\Bbb R^n\times \Bbb R^m)$ by iteration.
Similarly, we have $\big\|\sum_{j,k\in \Bbb Z} D^N_j\overline D^N_kM_bf\big\|_{L^2(\Bbb R^n\times \Bbb R^m)}\le
        C\|f\|_{L^2(\Bbb R^n\times \Bbb R^m)}$.

Since $\lim\limits_{j\to \infty} S_jM_{b_1} =I$,
      $\lim\limits_{j\to -\infty} S_jM_{b_1}=0$ on $L^2(\Bbb R^n)$ (cf. [DJS, p. 17]) and so
does $\overline S_k$ on $L^2(\Bbb R^m)$, we write
$$\allowdisplaybreaks\aligned
I&=\bigg(\sum_{j\in \Bbb Z} D_jM_{b_1}\bigg)\bigg(\sum_{j'\in \Bbb Z} D_{j'}M_{b_1}\bigg)\\
 &=\sum_{|\ell|>N}\sum_{j\in \Bbb Z}D_jM_{b_1}D_{j+\ell}M_{b_1}
       + \sum_{j\in \Bbb Z}\sum_{|\ell|\le N}D_jM_{b_1}D_{j+\ell}M_{b_1} \\
 &=\sum_{|\ell|>N}\sum_{j\in \Bbb Z}D_jM_{b_1}D_{j+\ell}M_{b_1} +  \sum_{j\in \Bbb Z}D_jM_{b_1}D^N_jM_{b_1} \\
 &:=\Cal R_N+\Cal V_N,
\endaligned                                                              \tag1.2$$
where $D^N_j:= \sum_{|\ell|\le N}D_{j+\ell}$. By [DJS, Lemma 2.2],
$\lim\limits_{N\to \infty}\Cal V_N=I$ in $L^2(\Bbb R^n)$ and
$\lim\limits_{N\to \infty}\|\Cal R_N\|_{L^2 \mapsto L^2}$$=0$, which
guarantees the existence of $\Cal V_N^{-1}$. For the product space,
we write
$$\allowdisplaybreaks\align
I&=\bigg(\sum_{|\ell_1|>N}\sum_{j\in \Bbb Z}D_jM_{b_1}D_{j+\ell_1}M_{b_1}
             +  \sum_{j\in \Bbb Z}D_jM_{b_1}D^N_jM_{b_1}\bigg) \\
 &\qquad\qquad\times
        \bigg(\sum_{|\ell_2|>N}\sum_{k\in \Bbb Z}\overline D_kM_{b_2}\overline D_{k+\ell_2}M_{b_2}
                   +  \sum_{k\in \Bbb Z}\overline D_kM_{b_2}\overline D^N_kM_{b_2}\bigg)\\
 &=\sum_{|\ell_1|>N}\sum_{j\in \Bbb Z}D_jM_{b_1}D_{j+\ell_1}M_{b_1}
      \sum_{|\ell_2|>N}\sum_{k\in \Bbb Z}\overline D_kM_{b_2}\overline D_{k+\ell_2}M_{b_2} \\
 &\qquad       + \sum_{|\ell_1|>N}\sum_{j\in \Bbb Z}D_jM_{b_1}D_{j+\ell_1}M_{b_1}
         \sum_{k\in \Bbb Z}\overline D_kM_{b_2}\overline D^N_kM_{b_2} \\
 &\qquad   + \sum_{j\in \Bbb Z}D_jM_{b_1}D^N_jM_{b_1}
       \sum_{|\ell_2|>N}\sum_{k\in \Bbb Z}\overline D_kM_{b_2}\overline D_{k+\ell_2}M_{b_2} \\
 &\qquad +\sum_{j\in \Bbb Z}D_jM_{b_1}D^N_jM_{b_1}\sum_{k\in \Bbb Z}\overline D_kM_{b_2}\overline D^N_kM_{b_2}, \\
 &:=R^1_N+R^2_N+R^3_N+V_N.
\endalign$$
By iteration, $V_N$ converges strongly on $L^2$ and $V_N^{-1}$ is bounded on $L^2$.
To use the $L^2$ boundedness of $V_N$ to get the $L^2$
boundedness of $T,$ we need to show that how $V_N$ does act on the test function. For this purpose,
let $\Lambda^s(\Bbb R^n\times \Bbb R^m)$ denote the closure of $C^\eta_0(\Bbb R^n\times \Bbb R^m)$ with respect to the norm
   $\|\cdot\|_s, 0<s<\eta$. The following lemma shows the properties of operators $V_N$ acting on $\Lambda^s$.

\proclaim{Lemma 3}
Let $b_1, b_2$ be para-accretive functions on $\Bbb R^n, \Bbb R^m$ respectively
and $\varepsilon$ be the common regularity exponent of the approximations to the identity associated to $b_1, b_2$.
Suppose $0<s<\varepsilon/2.$ Then
\roster
\item"(i)" $V_N=\sum_{j\in \Bbb Z}D_jM_{b_1}D^N_jM_{b_1}\sum_{k\in \Bbb Z}
             \overline D_kM_{b_2}\overline D^N_kM_{b_2}$ converges strongly on $\Lambda^s$;
\item"(ii)" $V_N$ is bounded on $\Lambda^s$;
\item"(iii)" $\|V_N-I\|_s\to 0$ as $N\to \infty$.
\endroster
\endproclaim
The proof of Lemma 3 will be given in section 2.

To see how one can use Lemma 3 to show Theorem 1, let
$b(x_1,x_2)=b_1(x_1)b_2(x_2),$ where $b_1$ and $b_2$ are para-accretive
functions on $\Bbb R^n$ and $\Bbb R^m,$ respectively, $f\in
\Lambda^s \cap L^2$, and $g\in C^s_0$. Suppose $\theta\in C^s_0$ and
$\theta=1$ on a open set containing the support of $g$, then we can
define $\langle bTbf, g \rangle =\langle bTb(\theta f), g\rangle
+\langle  bTb((1-\theta) f), g \rangle $. The first term makes sense
since $\theta f\in C^s_0$ and the second term can be defined by the
conditions of $(A_1)$ and $(A_2)$ and the assumption $f\in L^2.$
Hence $T$ can be extended to a continuous linear operator from
$b\Lambda^s \cap L^2$ into $(bC^s_0)^\prime.$ Assume Lemma 3 for the
moment,  $V_N$ is defined and bounded on $\Lambda^s\cap L^2$ with
the norm $\|\cdot\|_s+\|\cdot\|_{L^2}$, and it is invertible on
$\Lambda^s\cap L^2$ if $N$ large enough. Notice that $V_N$ converges
strongly on $L^2$ since
$\sup\limits_{L_1,L_2}\Big\|\sum\limits_{j,k=L_1}^{L_2}
D_jM_{b_1}D^N_jM_{b_1}
      \overline D_kM_{b_2}\overline D^N_kM_{b_2}\Big\|_{L^2\mapsto L^2}<\infty$
          and $V_Nf$ converges in $L^2$ if $f\in C^\eta_0$. Thus
$V_N$ converges strongly on $\Lambda^s\cap L^2$, by Lemma 3.
It is clear that $\Lambda^s\cap L^2$ is dense in $L^2$.
To prove Theorem 1, it suffices to show that
$$\big|\langle bg_0, Tbf_0 \rangle\big|\le C\|g_0\|_{L^2}\|f_0\|_{L^2}\qquad\text{for any}\ g_0, f_0\in \Lambda^s\cap L^2$$
Let $g_0\in \Lambda^s\cap L^2$ and let $g_1=V_N^{-1}g_0$ and set
$$U_{L_1,L_2}=\sum_{L_1\leq j\leq L_2}D_jM_{b_1}D^N_jM_{b_1}
        \sum_{L_1\leq k\leq L_2}\overline D_{k}M_{b_2}\overline D^N_{k}M_{b_2}.$$
Then $g_1\in \Lambda^s\cap L^2$ and
$\lim_{\Sb L_1\to -\infty\\ L_2\to +\infty\endSb}  U_{L_1,L_2}g_1=g_0$ in $\Lambda^s\cap L^2$.
Hence
$$\langle g_0, bTbf_0 \rangle =\lim\limits_{\Sb L_1\to -\infty \\ L_2\to +\infty \endSb}
                                                 \langle U_{L_1,L_2}g_1, bTbf_0 \rangle.$$
Similarly, let $f_0\in \Lambda^s\cap L^2$ and let $f_1=V_N^{-1}f_0$. Then $f_1\in \Lambda^s\cap L^2$
  and
$\lim_{\Sb L'_1\to -\infty\\ L'_2\to +\infty\endSb} U_{L'_1,L'_2}f_1=f_0$ in $\Lambda^s\cap L^2$.
Thus
$$\langle g_0, bTbf_0 \rangle =\lim\limits_{\Sb L_1\to -\infty \\ L_2\to +\infty \endSb}
               \lim\limits_{\Sb L'_1\to -\infty \\ L'_2\to +\infty \endSb}
                         \langle U_{L_1,L_2}g_1, bTbU_{L'_1,L'_2}f_1 \rangle.$$
Therefore, we have to show that $U_{L_1,L_2}bTbU_{L'_1,L'_2}$ is bounded on $L^2$ as $L_1,L'_1\to -\infty$ and $L_2,L_2'\to \infty$.
We remark that the $L^2$ boundedness of $U_{L_1,L_2}bTbU_{L'_1,L'_2}$ uniformly for
$L_1,L_2,L^\prime_1,L^\prime_2$ was proved in [DJS] under the assumptions
that $T(1)={}^tT(1)=\widetilde T(1)={}^t\widetilde T(1)=0.$ However, we will show this result without assuming
$T(1)={}^tT(1)=\widetilde T(1)={}^t\widetilde T(1)=0$ and the proof will be given in section 3.
Same as the relation between $S_j$ and $S_j(x,y)$, if $D$ is an operator,
then we use $D(x,y)$ denotes its corresponding kernel through the article,
and the same remark apply to $D^N_j$, $P_j$, and so on.
For simplicity, we also denote $\int dv$ by $\int_{\Bbb R^n\times\Bbb R^m} dv_1dv_2$
and similarly for other variables.

\head \S2. Proof of Lemma 3
\endhead
For $f$ defined on $\Bbb R^n\times \Bbb R^m$, we use $\|f\|_{(L^\infty(\Bbb R^n), \lambda^\beta(\Bbb R^m))}$
and $\|f\|_{(\lambda^\beta(\Bbb R^n), L^\infty(\Bbb R^m))}$
to express
$$\sup\limits_{\Sb x_1\in \Bbb R^n\\ x_2\ne y_2 \endSb} \frac {|f(x_1,x_2)-f(x_1,y_2)|}{|x_2-y_2|^\beta}
\qquad\text{and}\qquad \sup\limits_{\Sb x_1\ne y_1\\ x_2\in \Bbb R^m\endSb}
                              \frac {|f(x_1,x_2)-f(y_1,x_2)|}{|x_1-y_1|^\beta},$$
respectively. To prove Lemma 3, we need the following estimates for
$\{D_jM_{b_1}\overline D_{k}M_{b_2}\}$.
\proclaim{Lemma 4} Let $b_1, b_2$ be para-accretive functions on $\Bbb R^n, \Bbb R^m$ and
$\varepsilon$ be the regularity exponent of the approximations to
the identity associated to $b_1, b_2$ respectively.
For $f\in \Lambda_s(\Bbb R^n\times \Bbb R^m), s<\varepsilon$,
\roster
\item"(i)" $\|D_jM_{b_1}\overline D_{k}M_{b_2}f\|_{L^\infty} \le C2^{-(j+k)s}\|f\|_s$,
\item"(ii)"  $\|D_jM_{b_1}\overline D_{k}M_{b_2}f\|_{(L^\infty(\Bbb R^n), \lambda^\beta(\Bbb R^m))}
               \le C2^{-js}2^{k(\beta-s)}\|f\|_s$\qquad if\ \ $0<s\leq\beta< \varepsilon$,
\item"(iii)"  $\|D_jM_{b_1}\overline D_{k}M_{b_2}f\|_\beta\le C2^{(j+k)(\beta-s)}\|f\|_s$
                 \qquad if\ \ $0<s\leq\beta< \varepsilon$,
\item"(iv)" $\|D^N_jM_{b_1}\overline D^N_{k}M_{b_2}f\|_s\le CN^2\|f\|_s$.
\endroster
\endproclaim
\demo{Proof} For (i), the cancellations of $D_j$ and $\overline D_k$
give
$$\align &D_jM_{b_1}\overline D_{k}M_{b_2}f(x_1,x_2)\\
&=\int D_j(x_1,y_1)b_1(y_1)\overline D_{k}(x_2,y_2)b_2(y_2)
     \big[f(y_1,y_2)-f(x_1,y_2)-f(y_1,x_2)+f(x_1.x_2)\big]dy.
\endalign$$
Since $S_j(x,y)=0$ for $|x-y|\ge c2^{-j}$, the size conditions of
$D_j$ and $\overline D_k$ yield
$$\align |D_jM_{b_1}\overline D_{k}M_{b_2}f(x_1,x_2)|
         &\le C\|f\|_s \int_{\Sb |x_1-y_1|< c2^{-j} \\ |x_2-y_2|<c2^{-k} \endSb}
                     2^{jn+km}|x_1-y_1|^s|x_2-y_2|^sdy_1dy_2 \\
         &\le C2^{-(j+k)s}\|f\|_s.
\endalign$$

To obtain (ii), we write
$$\align
&D_jM_{b_1}\overline D_{k}M_{b_2}f(x_1,x_2)-D_jM_{b_1}\overline D_{k}M_{b_2}f(x_1,y_2) \\
&\qquad= \int D_j(x_1,z_1)b_1(z_1)\big[\overline D_{k}(x_2,z_2)-\overline D_{k}(y_2,z_2)\big]b_2(z_2)f(z_1,z_2)dz.
\endalign$$
If $|x_2-y_2|\le c2^{-k}$,
the cancellations of $D_j$ and $\overline D_{k}$ yield that
$$\align
&|D_jM_{b_1}\overline D_{k}M_{b_2}f(x_1,x_2)-D_jM_{b_1}\overline D_{k}M_{b_2}f(x_1,y_2)| \\
&\qquad= \Big|\int D_j(x_1,z_1)b_1(z_1)\big[\overline D_{k}(x_2,z_2)-\overline D_{k}(y_2,z_2)\big]b_2(z_2) \\
&\hskip3cm \times            \big[f(z_1,z_2)-f(x_1,z_2)-f(z_1,x_2)+f(x_1,x_2)\big]dz\Big|\\
&\qquad\le C\bigg(\int_{\Sb |x_1-z_1|< c2^{-j} \\ |x_2-z_2|<c2^{-k} \endSb}
   +\int_{\Sb |x_1-z_1|< c2^{-j} \\ |y_2-z_2|<c2^{-k} \endSb}\bigg)
                 2^{jn}|x_2-y_2|^\beta2^{k(m+\beta)}|x_1-z_1|^s |x_2-z_2|^s \|f\|_s dz_1dz_2.\\
\endalign$$
Hence, for $|x_2-y_2|\le c2^{-k}$,
$$|D_jM_{b_1}\overline D_{k}M_{b_2}f(x_1,x_2)-D_jM_{b_1}\overline D_{k}M_{b_2}f(x_1,y_2)|
              \le C|x_2-y_2|^\beta2^{-js}2^{k(\beta-s)}\|f\|_s. \tag2.1$$
For $|x_2-y_2|> c2^{-k}$, (i) gives
$$|D_jM_{b_1}\overline D_{k}M_{b_2}f(x_1,x_2)-D_jM_{b_1}\overline D_{k}M_{b_2}f(x_1,y_2)|
   \le C2^{-(j+k)s}\|f\|_s\le C|x_2-y_2|^\beta2^{-js}2^{k(\beta-s)}\|f\|_s.$$
Therefore, we obtain (ii).

To estimate (iii), we write
$$\align
&D_kM_{b_1}\overline D_{k'}M_{b_2}f(x_1,x_2)-D_kM_{b_1}\overline D_{k'}M_{b_2}f(x_1,y_2) \\
&\qquad\qquad      -D_kM_{b_1}\overline D_{k'}M_{b_2}f(y_1,x_2)
                             +D_kM_{b_1}\overline D_{k'}M_{b_2}f(y_1,y_2) \\
&\qquad= \int\big[D_k(x_1,z_1)-D_k(y_1,z_1)\big]b_1(z_1)
             \big[\overline D_{k'}(x_2,z_2)-\overline D_{k'}(y_2,z_2)\big]b_2(z_2)f(z_1,z_2)dz.
\endalign$$
For $|x_1-y_1|\le c2^{-j}$ and $|x_2-y_2|\le c2^{-k}$, we use
the the cancellations of $D_k$ and $\overline D_{k'}$ to get
$$\align
&\big|D_jM_{b_1}\overline D_{k}M_{b_2}f(x_1,x_2)-D_jM_{b_1}\overline D_{k}M_{b_2}f(x_1,y_2) \\
&\qquad\qquad      -D_jM_{b_1}\overline D_{k}M_{b_2}f(y_1,x_2)+D_jM_{b_1}\overline D_{k}M_{b_2}f(y_1,y_2)\big| \\
&\qquad=\bigg| \int \big[D_j(x_1,z_1)-D_j(y_1,z_1)\big]b_1(z_1)\big[\overline D_{k}(x_2,z_2)-\overline D_{k}(y_2,z_2)\big]b_2(z_2)\\
&\hskip3cm \times \big[f(z_1,z_2)-f(x_1,z_2)-f(z_1,x_2)+f(x_1,x_2)\big]dz\bigg|\\
&\qquad\le C\bigg(\int_{\Sb |x_1-z_1|< c2^{-j} \\ |x_2-z_2|<c2^{-k}
\endSb}+\int_{\Sb |x_1-z_1|< c2^{-j} \\ |y_2-z_2|<c2^{-k} \endSb}
                +\int_{\Sb |y_1-z_1|< c2^{-j} \\ |x_2-z_2|<c2^{-k} \endSb}
                +\int_{\Sb |y_1-z_1|< c2^{-j} \\ |y_2-z_2|<c2^{-k} \endSb}\bigg)\\
&\hskip3cm\times
|x_1-y_1|^\beta2^{j(n+\beta)}|x_2-y_2|^\beta2^{k(m+\beta)}|x_1-z_1|^s
                                |x_2-z_2|^s \|f\|_s dz_1dz_2.\\
&\qquad\le C|x_1-y_1|^\beta|x_2-y_2|^\beta
2^{(j+k)(\beta-s)}\|f\|_s.
\endalign$$
For $|x_1-y_1|\le c2^{-j}$ and $|x_2-y_2|> c2^{-k}$, (2.1) shows that
$$\align
&\big|D_jM_{b_1}\overline D_{k}M_{b_2}f(x_1,x_2)-D_jM_{b_1}\overline D_{k}M_{b_2}f(x_1,y_2) \\
&\qquad\qquad      -D_jM_{b_1}\overline D_{k}M_{b_2}f(y_1,x_2)+D_jM_{b_1}\overline D_{k}M_{b_2}f(y_1,y_2)\big| \\
&\qquad\le \big|D_jM_{b_1}\overline D_{k}M_{b_2}f(x_1,x_2)-D_jM_{b_1}\overline D_{k}M_{b_2}f(y_1,x_2)\big| \\
&\qquad\qquad      +\big|D_jM_{b_1}\overline D_{k}M_{b_2}f(x_1,y_2)-D_jM_{b_1}\overline D_{k}M_{b_2}f(y_1,y_2)\big| \\
&\qquad\le C|x_1-y_1|^\beta2^{j(\beta-s)}2^{-ks}\|f\|_s\\
&\qquad\le
C|x_1-y_1|^\beta|x_2-y_2|^\beta2^{j(\beta-s)}2^{k(\beta-s)}\|f\|_s.
\endalign$$
The case $|x_1-y_1|> c2^{-j}$ and $|x_2-y_2|\le c2^{-k}$ is similar, so we consider the final case
  $|x_1-y_1|> c2^{-j}$ and $|x_2-y_2|> c2^{-k}$. (i) gives
$$\align
&\big|D_jM_{b_1}\overline D_{k}M_{b_2}f(x_1,x_2)-D_jM_{b_1}\overline D_{k}M_{b_2}f(x_1,y_2) \\
&\qquad\qquad      -D_jM_{b_1}\overline D_{k}M_{b_2}f(y_1,x_2)+D_jM_{b_1}\overline D_{k}M_{b_2}f(y_1,y_2)\big| \\
&\qquad\le C2^{-(j+k)s}\|f\|_s \\
&\qquad\le C|x_1-y_1|^\beta|x_2-y_2|^\beta
2^{(j+k)(\beta-s)}\|f\|_s.
\endalign$$
The estimate of (iii) is completed.

The estimate of (iv) can be done
by the same argument as (iii).$\hfill$\qed
\enddemo

We now return to show Lemma 3.

\demo{Proof of Lemma 3} Suppose that $f\in \Lambda_s(\Bbb R^n\times \Bbb R^m).$
Set $G_{jk}f:=D_jM_{b_1}\overline
D_{k}M_{b_2}D^N_jM_{b_1}\overline D^N_{k}M_{b_2}f$. Given
$(x_1,x_2), (y_1,y_2)\in \Bbb R^n\times\Bbb R^m$, choose $j_0,
k_0\in \Bbb Z$ such that
  $2^{-j_0}\le |x_1-y_1|\le 2^{-j_0+1}$ and $2^{-k_0}\le |x_2-y_2|\le 2^{-k_0+1}$. Then
Lemma 4 implies that
$$\aligned
&\bigg|\sum_{j,k} \big[G_{jk}f(x_1,x_2)-G_{jk}f(y_1,x_2)-G_{jk}f(x_1,y_2)+G_{jk}f(y_1,y_2)\big]\bigg|\\
&\qquad\le \bigg(\sum_{j\ge j_0, k\ge k_0} + \sum_{j\ge j_0, k< k_0} +\sum_{j< j_0, k\ge k_0} +\sum_{j< j_0, k< k_0}\bigg) \\
&\qquad\qquad |G_{jk}f(x_1,x_2)-G_{jk}f(y_1,x_2)-G_{jk}f(x_1,y_2)+G_{jk}f(y_1,y_2)|\\
&\qquad\le \sum_{j\ge j_0, k\ge k_0} 4\|G_{jk}f\|_{L^\infty}
              + \sum_{j\ge j_0, k< k_0} 2|x_2-y_2|^\beta\|G_{jk}f\|_{(L^\infty(\Bbb R^n), \lambda^\beta(\Bbb R^m))}\\
&\qquad\quad+  \sum_{j< j_0, k\ge k_0} 2|x_1-y_1|^\beta\|G_{jk}f\|_{(\lambda^\beta(\Bbb R^n), L^\infty(\Bbb R^m))} \\
&\qquad \quad            + \sum_{j< j_0, k< k_0} |x_1-y_1|^\beta|x_2-y_2|^\beta\|G_{kk'}f\|_\beta\\
&\qquad \le C_N2^{-(j_0+k_0)s}\|f\|_s
+C_N2^{-j_0s}2^{k_0(\beta-s)}|x_2-y_2|^\beta\|f\|_s    \\
&\qquad \quad +C_N2^{j_0(\beta-s)} |x_1-y_1|^\beta2^{-k_0s}
              +C_N2^{(j_0+k_0)(\beta-s)}|x_1-y_1|^\beta|x_2-y_2|^\beta\|f\|_s \\
&\qquad\le C_N|x_1-y_1|^s|x_2-y_2|^s\|f\|_s.
\endaligned \tag2.2 $$
Hence $\big\|\sum_{j,k} D_jM_{b_1}\overline D_{k}M_{b_2}D^N_jM_{b_1}\overline D^N_{k}M_{b_2}f\big\|_s
              \le C_N\|f\|_s$. If $f\in C^\eta_0$ for some $\eta>s$,
then the series $\sum_{j,k}D_jM_{b_1}\overline D_{k}M_{b_2}D^N_jM_{b_1}\overline D^N_{k}M_{b_2}f$ converges
uniformly and in $\Lambda^s$ norm. This implies $V_N$ is bounded on $\Lambda^s$.
We now show that $\|R_N^i\|_s\to 0$ as $N\to \infty$ for $i=1,2,3.$
We only show this limit for $i=2$ because the proofs for $i=1,3$ are similar and we leave the details to the readers.
We first rewrite
$$\align
R_N^2f
&=\sum_{|\ell_1|>N}\sum_{j\in \Bbb Z}D_jM_{b_1}D_{j+\ell_1}M_{b_1}
         \sum_{k\in \Bbb Z}\overline D_kM_{b_2}\overline D^N_kM_{b_2}f \\
&=\bigg(\sum_{j\in \Bbb Z}D_jM_{b_1}(I-S_{j+N}M_{b_1})
             +\sum_{j\in \Bbb Z}D_jM_{b_1}S_{j-N-1}M_{b_1}
   \bigg)\sum_{k\in \Bbb Z}\overline D_kM_{b_2}\overline D^N_kM_{b_2}f\\
&=\sum_{j,k\in \Bbb Z}D_jM_{b_1}(I-S_{j+N}M_{b_1})\overline D_kM_{b_2}\overline D^N_kM_{b_2}f
       +\sum_{j,k\in \Bbb Z}D_jM_{b_1}S_{j-N-1}M_{b_1}\overline D_kM_{b_2}\overline D^N_kM_{b_2}f
\endalign$$
Since $\int_{\Bbb R^n} S_j(x_1,y_1)b_1(y_1)dy_1=1, j\in \Bbb Z$, we have
$$\align &(I-S_{j+N}M_{b_1})\overline D^N_kM_{b_2}f(x_1,x_2)\\
&\qquad=\overline D^N_kM_{b_2}f(x_1,x_2)-S_{j+N}M_{b_1}\overline D^N_kM_{b_2}f(x_1,x_2)\\
&\qquad=\int_{\Bbb R^m} D^N_k(x_2,z_2)b_2(z_2)f(x_1,z_2)dz_2
                 -\int S_{j+N}(x_1,z_1)b_1(z_1)\overline D^N_k(x_2,z_2)b_2(z_2)f(z_1,z_2)dz\\
&\qquad=\int S_{j+N}(x_1,z_1)b_1(z_1)\overline D^N_k(x_2,z_2)b_2(z_2)\big[f(x_1,z_2)-f(z_1,z_2)\big]dz.
\endalign$$
Hence we can regard $(I-S_{j+N}M_{b_1})$ as $D_{j+N}M_{b_1}$ so that
$$\|(I-S_{j+N}M_{b_1})\overline D^N_kM_{b_2}\|_{(L^\infty(\Bbb R^n), \lambda^s(\Bbb R^m))}
       \le CN2^{-(N+j)s}\|f\|_s.\tag 2.3$$
By the same argument of Lemma 4, we obtain
$$\align
&\|D_jM_{b_1}\overline D_kM_{b_2}f\|_{L^\infty}\le C2^{-ks}\|f\|_{(L^\infty(\Bbb R^n), \lambda^s(\Bbb R^m))}\\
&\|D_jM_{b_1}\overline D_kM_{b_2}f\|_{(L^\infty(\Bbb R^n), \lambda^\beta(\Bbb R^m))}
           \le C2^{k(\beta-s)}\|f\|_{(L^\infty(\Bbb R^n), \lambda^s(\Bbb R^m))} \\
&\|D_jM_{b_1}\overline D_kM_{b_2}f\|_{(\lambda^\beta(\Bbb R^n), L^\infty(\Bbb R^m))}
           \le C2^{j\beta}2^{-ks}\|f\|_{(L^\infty(\Bbb R^n), \lambda^s(\Bbb R^m))} \\
&\|D_jM_{b_1}\overline D_kM_{b_2}f\|_\beta
           \le C2^{j\beta}2^{k(\beta-s)}\|f\|_{(L^\infty(\Bbb R^n), \lambda^s(\Bbb R^m))}.
\endalign$$
The above estimates together with (2.3) and the same method of (2.2) show that
$$\bigg\|\sum_{j,k\in \Bbb Z}D_jM_{b_1}(I-S_{j+N}M_{b_1})\overline D_kM_{b_2}\overline D^N_kM_{b_2}\bigg\|_s
  \le CN2^{-Ns}\|f\|_s$$
and hence $N2^{-N}\to 0$ as $N\to \infty$.
Let $H^N_{j,k}f=D_jM_{b_1}S_{j-N-1}M_{b_1}\overline D_kM_{b_2}\overline D^N_kM_{b_2}f$
  and let $H^N_j(x_1,y_1)$ be the kernel of $D_jM_{b_1}S_{j-N-1}M_{b_1}$.
Then $\int_{\Bbb R^n} H^N_j(x_1,y_1)dy_1=D_jM_{b_1}S_{j-N-1}M_{b_1}(1)$ $=D_jM_{b_1}(1)=0$
and $H^N_j(x_1,y_1)=0$ if $|x_1-y_1|\ge C2^{-(j-N)}$.
By the cancellation of $D_j$ and the estimate of $S_{j-N-1}$,
$$\aligned
|H^N_{j}(x_1,y_1)|
&=\bigg|\int
     D_j(x_1,z_1)b_1(z_1)\big[S_{j-N-1}(z_1,y_1)-S_{j-N-1}(x_1,y_1)\big]b_1(y_1)dz_1\bigg|\\
&\le C\int_{|x_1-z_1|\le c2^{-j}}
        2^{jn}|x_1-z_1|^\varepsilon 2^{(j-N-1)(n+\varepsilon)}dz_1\\
&\le C2^{-N\varepsilon}2^{(j-N)n}.
\endaligned$$
Put $\overline D^N_kM_{b_2}f=h$. The above inequality and the cancellations yield that
$$\align
&|H^N_{j,k}f(x_1,x_2)|\\
&\quad=\bigg|\int H_j^N(x_1,y_1)\overline D_k(x_2,y_2)b_2(y_2)h(y_1,y_2)dy\bigg|\\
&\quad=\bigg|\int H_j^N(x_1,y_1)\overline D_k(x_2,y_2)b_2(y_2)\big[h(y_1,y_2)-h(x_1,y_2)-h(y_1,x_2)+h(x_1,x_2)\big]dy\bigg| \\
&\quad\le C\int_{\Sb |x_1-y_1|< C2^{-(j-N)} \\ |x_2-y_2|<c2^{-k} \endSb}
                      | H_j^N(x_1,y_1)|2^{km}|x_1-y_1|^s|x_2-y_2|^s\|h\|_sdy_1dy_2 \\
&\quad\le C2^{-N\varepsilon}2^{-(j-N)s}2^{-ks}\|h\|_s \\
&\quad\le CN2^{-N\varepsilon}2^{-(j-N)s}2^{-ks}\|f\|_s
\endalign$$
and then
$$\|H^N_{j,k}f\|_{L^\infty}\le CN2^{-N\varepsilon}2^{-(j-N)s}2^{-ks}\|f\|_s. \tag 2.4$$
If $|x_1-x_1'|\le C2^{-(j-N)}$, then the size conditions show that
$$\allowdisplaybreaks\aligned
|H^N_j(x_1,y_1)-H^N_j(x'_1,y_1)|
&= \bigg|\int_{\Bbb R^n} \big[D_j(x_1,z_1)-D_j(x'_1,z_1)\big]S_{j-N-1}(z_1,y_1)dz_1\bigg| \\
&\le C \int_{\Sb {\{z_1\in \Bbb R^n: |x_1-z_1|\le c2^{-j}}\\  {\text{or}\ |x'_1-z_1|\le c2^{-j}\}}\endSb}
        |x_1-x'_1|^\varepsilon 2^{j(n+\varepsilon)}2^{(j-N)n}dz_1 \\
&\le C2^{j\varepsilon}|x_1-x'_1|^\varepsilon 2^{(j-N)n}.
\endaligned \tag2.5$$
For $|x_1-y_1|\le C2^{-(j-N)}$, we use (2.5) to get
$$\allowdisplaybreaks\align
&|H^N_{j,k}f(x_1,x_2)-H^N_{j,k}f(y_1,x_2)| \\
&\qquad=\bigg|\int \big[H^N_j(x_1,z_1)-H^N_j(y_1,z_1)\big]\overline D_k(x_2,z_2)b_2(z_2)h(z_1,z_2)dz\bigg|\\
&\qquad=\bigg|\int \big[H^N_j(x_1,z_1)-H^N_j(y_1,z_1)\big]\overline D_k(x_2,z_2)b_2(z_2) \\
&\qquad\qquad\times \big[h(z_1,z_2)-h(x_1,z_2)-h(z_1,x_2)+h(x_1,x_2)\big]dz\bigg|\\
&\qquad \le C\bigg(\int_{\Sb |x_1-z_1|\le C2^{-(j-N)} \\ |x_2-z_2|\le c2^{-k} \endSb}
                   + \int_{\Sb |y_1-z_1|\le C2^{-(j-N)} \\ |x_2-z_2|\le c2^{-k} \endSb} \bigg)\\
&\hskip3cm\times     2^{j\varepsilon}|x_1-y_1|^\varepsilon 2^{(j-N)n}2^{km}|x_1-z_1|^s|x_2-z_2|^sdz_1dz_2\\
&\qquad\le C2^{j\varepsilon}2^{-(j-N)s}|x_1-y_1|^\varepsilon2^{-ks}\|h\|_s\\
&\qquad\le CN2^{j\varepsilon}2^{-(j-N)s}|x_1-y_1|^\varepsilon2^{-ks}\|f\|_s,
\endalign$$
where $h=\overline D^N_kM_{b_2}f$.
For $|x_1-y_1|> C2^{-(j-N)}$, (2.4) implies
$$|H^N_j(x_1,y_1)-H^N_j(x'_1,y_1)| \le CN2^{-N\varepsilon}2^{-(j-N)s}2^{-ks}\|f\|_s\le
  CN2^{j\varepsilon}2^{-(j-N)s}|x_1-y_1|^\varepsilon2^{-ks}\|f\|_s$$
so that
$$\|H^N_{j,k}f\|_{(\lambda^\varepsilon(\Bbb R^n), L^\infty(\Bbb R^m))}
                          \le CN2^{j\varepsilon}2^{-(j-N)s}2^{-ks}\|f\|_s.\tag2.6$$
Using the fact $\|f\|_{(\lambda^\varepsilon(\Bbb R^n), L^\infty(\Bbb R^m))}
              \le \|f\|_{L^\infty}^{\frac {\varepsilon-\beta}\varepsilon}
                  \|f\|_{(\lambda^\varepsilon(\Bbb R^n), L^\infty(\Bbb R^m))}^{\frac \beta\varepsilon}$,
$0<\beta<\varepsilon$, (2.4) and (2.6) give
$$\|H^N_{j,k}f\|_{(\lambda^\beta(\Bbb R^n), L^\infty(\Bbb R^m))}
                   \le CN2^{-N(\varepsilon-2\beta)}2^{(j-N)(\beta-s)}2^{-ks}\|f\|_s.\tag2.7$$
Similarly, we have
$$\|H^N_{j,k}f\|_{(L^\infty(\Bbb R^n, \lambda^\beta(\Bbb R^m)))}
                           \le CN2^{-N\varepsilon}2^{-(j-N)s}2^{k(\beta-s)}\|f\|_s \tag2.8$$
and
$$\|H^N_{j,k}f\|_\beta \le  CN2^{-N(\varepsilon-2\beta)}2^{(j-N)(\beta-s)}2^{k(\beta-s)}\|f\|_s.\tag2.9$$
Plugging $\tilde H^N_{j,k}f=2^{N(\varepsilon-2\beta)}H^N_{j+N,k}f$ in (2.4) and (2.7)-(2.9),
we use the same method of (2.2) to have
$$\Big\|\sum_{j,k}H^N_{j,k}f\Big\|_s\le CN2^{-N(\varepsilon-2\beta)}\|f\|_s\qquad \text{for
$s<\beta<\varepsilon$}.$$ If $s<\varepsilon/2$, we can choose $\beta$
so that $N2^{-N(\varepsilon-2\beta)}\to 0$ as $N\to \infty.$ The proof is finished. $\hfill$\qed
\enddemo

\head
\S3. Proof of Theorem 1
\endhead

To prove Theorem 1, it suffices to show
$$\lim\limits_{\Sb L_1\to -\infty \\ L_2\to +\infty \endSb}\langle U_{L_1,L_2}g, bTM_bU_{L'_1,L'_2}f \rangle \le C\|f\|_{L^2}\|g\|_{L^2},$$
where $C$ is independent of $f,g\in \Lambda^s\cap L^2$.
For simplicity, We use $D_{jk}M_b$ and $D^N_{jk}M_b$ to express
$D_jM_{b_1}\overline D_kM_{b_2}$ and $D^N_jM_{b_1}\overline D^N_kM_{b_2}$, respectively,
and recall that $\int dv$  denotes $\int_{\Bbb R^n\times\Bbb R^m} dv_1dv_2$.
For $f,g\in \Lambda^s\cap L^2$,
$$\align
&\lim\limits_{\Sb L_1\to -\infty \\ L_2\to +\infty \endSb}\langle U_{L_1,L_2}g, bTM_bU_{L'_1,L'_2}f \rangle \\
&\quad=\sum_{j,k,j',k'}\Big\langle M_bD^N_{jk}M_b g,
          {}^t\!D_jM_{b_1}{}^t\overline
        D_kM_{b_2}TM_bD_{j'k'}M_bD^N_{j'k'}M_bf\Big\rangle \\
&\quad= \sum_{j,k}
      \bigg(\sum_{\Sb j'\ge j \\ k'\ge k \endSb}+\sum_{\Sb j'\ge j \\ k'< k\endSb}
       +\sum_{\Sb j'< j \\ k'\ge k \endSb}+\sum_{\Sb j'< j\\ k'< k\endSb}\bigg)
        \bigg\langle M_bD^N_{jk}M_b g, {}^t\!D_jM_{b_1}{}^t\overline
        D_kM_{b_2}TM_bD_{j'k'}M_bD^N_{j'k'}M_bf\bigg\rangle \\
&\quad:=J_1+J_2+J_3+J_4.
\endalign$$

Since $J_3$ and $J_4$ are symmetric with respect to $J_2$ and $J_1$, respectively, we only prove the cases
$J_1$ and $J_2$. Using the one parameter idea, we directly minus
functions which we want such that the almost orthogonality argument
holds. Hence we have three remainder terms to handle. The first two
terms are mixed terms, we have the almost orthogonality argument for
one parameter but not two parameters. We do the estimates by the
following Theorem 5. The final remainder term looks like the
paraproduct $\Pi_b$ but easier. The terms do not satisfy the
conditions $\Pi_b(1)=b$ and $\Pi^*_b(1)=0.$ To be precise, we
write down $J_1$ as follows.

$$\allowdisplaybreaks\align
J_1 &=\sum_{j,k}\sum_{\Sb j'\ge j \\ k'\ge k \endSb} \bigg\langle
M_bD^N_{jk}M_b g, \int {}^t\!D_jM_{b_1}{}^t\overline
D_kM_{b_2}TM_bD_{j'}(\cdot,v_1)\overline D_{k'}(\cdot,v_2)M_bD^N_{j'k'}M_bf(v_1,v_2)dv\bigg\rangle \\
&=\sum_{j,k}\sum_{\Sb j'\ge j \\ k'\ge k \endSb} \bigg\langle
M_bD^N_{jk}M_b g, \iint
\big[{}^t\!D_j(\cdot,z_1)-{}^t\!D_j(\cdot,v_1)\big]b_1(z_1)
        \big[{}^t\overline D_k(\cdot,z_2)-{}^t\overline D_k(\cdot,v_2)\big]b_2(z_2)\\
&\hskip3cm\times TM_bD_{j'}(z_1,v_1)\overline D_{k'}
                     (z_2,v_2)M_bD^N_{j'k'}M_bf(v_1,v_2)dzdv\bigg\rangle \\
&\qquad+\sum_{j,k}\sum_{\Sb j'\ge j \\ k'\ge k \endSb}
          \bigg\langle M_bD^N_{jk}M_b g, \iint  {}^t\!D_j(\cdot,v_1)b_1(z_1)
        \big[{}^t\overline D_k(\cdot,z_2)-{}^t\overline D_k(\cdot,v_2)\big]b_2(z_2) \\
&\hskip3cm\times TM_bD_{j'}(z_1,v_1)\overline D_{k'}(z_2,v_2)M_bD^N_{j'k'}M_bf(v_1,v_2)dzdv\bigg\rangle  \\
&\qquad+\sum_{j,k}\sum_{\Sb j'\ge j \\ k'\ge k \endSb}
          \bigg\langle M_bD^N_{jk}M_b g, \iint \big[{}^t\!D_j(\cdot,z_1)-{}^t\!D_j(\cdot,v_1)\big]b_1(z_1)
        {}^t\overline D_k(\cdot,v_2)b_2(z_2) \\
&\hskip3cm\times TM_bD_{j'}(z_1,v_1)\overline D_{k'}(z_2,v_2)M_bD^N_{j'k'}M_bf(v_1,v_2) dzdv\bigg\rangle   \\
&\qquad+\sum_{j,k}\sum_{\Sb j'\ge j \\ k'\ge k \endSb}
          \bigg\langle M_bD^N_{jk}M_b g, \iint {}^t\!D_j(\cdot,v_1)b_1(z_1){}^t\overline D_k(\cdot,v_2)b_2(z_2)\\
&\hskip3cm\times TM_bD_{j'}(z_1,v_1)\overline D_{k'}(z_2,v_2)M_bD^N_{j'k'}M_bf(v_1,v_2)dzdv\bigg\rangle\\
&:=J_{11}+J_{12}+J_{13}+J_{14}.
\endalign$$

By an almost orthogonality argument, for $j'\ge j$ and $k'\ge k$,
$$\allowdisplaybreaks\align
&\bigg|\int \big[{}^t\!D_j(y_1,z_1)-{}^t\!D_j(y_1,v_1)\big]b_1(z_1)
        \big[{}^t\overline D_k(y_2,z_2)-{}^t\overline D_k(y_2,v_2)\big]b_2(z_2)\\
&\hskip5.5cm \times TM_bD_{j'}(z_1,v_1)\overline D_{k'}(z_2,v_2)dz\bigg| \\
&\qquad\le C2^{-(j'-j)\varepsilon'}2^{-(k'-k)\varepsilon'}\frac{2^{-j\varepsilon'}}{(2^{-j}+|y_1-v_1|)^{n+\varepsilon'}}
      \frac{2^{-k\varepsilon'}}{(2^{-k}+|y_2-v_2|)^{m+\varepsilon'}},
\endalign$$
where $0<\varepsilon'<\varepsilon$,
and hence
$$\allowdisplaybreaks\align
 |J_{11}|&\le C \sum_{j.k}\sum_{\Sb j'\ge j\\ k'\ge k \endSb}\iint
       2^{-(j'-j)\varepsilon'}2^{-(k'-k)\varepsilon'}
      \frac{2^{-j\varepsilon'}}{(2^{-j}+|y_1-v_1|)^{n+\varepsilon'}}
      \frac{2^{-k\varepsilon'}}{(2^{-k}+|y_2-v_2|)^{m+\varepsilon'}}\\
  &\hskip5cm\times |D^N_{jk}M_b g(y_1,y_2)|
      |D^N_{j'k'}M_bf(v_1,v_2)|dvdy.
 \endalign$$
Using Schwarz's inequality and the fact
$$\int_{\Bbb R^n}\frac{2^{-j\varepsilon'}}{(2^{-j}+|y_1-v_1|)^{n+\varepsilon'}}dy_1<\infty,$$
we obtain
$$\allowdisplaybreaks\align
&\int_{\Bbb R^n\times\Bbb R^n}  \frac{2^{-j\varepsilon'}}{(2^{-j}+|y_1-v_1|)^{n+\varepsilon'}}
    |D^N_{jk}M_b g(y_1,y_2)| |D^N_{j'k'}M_bf(v_1,v_2)|dy_1dv_1\\
&\le C\int_{\Bbb R^n} \bigg(\int_{\Bbb R^n} \frac{2^{-j\varepsilon'}}{(2^{-j}+|y_1-v_1|)^{n+\varepsilon'}}
            |D^N_{jk}M_b g(y_1,y_2)|^2dy_1\bigg)^{1/2}|D^N_{j'k'}M_bf(v_1,v_2)|dv_1 \\
&\le C\bigg(\int_{\Bbb R^n\times\Bbb R^n}\hskip-0.1cm \frac{2^{-j\varepsilon'}}{(2^{-j}+|y_1-v_1|)^{n+\varepsilon'}}
       |D^N_{jk}M_b g(y_1,y_2)|^2dy_1dv_1\bigg)^{1/2}
   \bigg(\int_{\Bbb R^n} |D^N_{j'k'}M_bf(v_1,v_2)|^2dv_1\bigg)^{1/2} \\
&\le C\bigg(\int_{\Bbb R^n}|D^N_{jk}M_b g(y_1,y_2)|^2dy_1\bigg)^{1/2}
         \bigg(\int_{\Bbb R^n} |D^N_{j'k'}M_bf(v_1,v_2)|^2dv_1\bigg)^{1/2}.
\endalign$$
Thus, $|J_{11}|$ is dominated by
$$\align &\int_{\Bbb R^m\times \Bbb R^m}\hskip-0.1cm \sum_k\sum_{k'\ge k} 2^{-(k'-k)\varepsilon'}
 \frac{2^{-k\varepsilon'}}{(2^{-k}+|y_2-v_2|)^{m+\varepsilon'}}
 \bigg(\int_{\Bbb R^n}\hskip-0.1cm \sum_j\sum_{j'\ge j}  2^{-(j'-j)\varepsilon'}
                          |D^N_{jk}M_bg(y_1,y_2)|^2dy_1\bigg)^{1/2} \\
 &\hskip3cm\times  \bigg(\int_{\Bbb R^n} \sum_j\sum_{j'\ge j} 2^{-(j'-j)\varepsilon'}
                       |D^N_{j'k'}M_bf(v_1,v_2)|^2dv_1\bigg)^{1/2}dy_2dv_2.
\endalign$$
We do a similar argument for variables $y_2$ and $v_2$ again to get
$$\align
|J_{11}|&\le C \bigg(\int \sum_{j,k} |D^N_{jk}M_b
g(y_1,y_2)|^2dy\bigg)^{1/2}
            \bigg(\int \sum_{j',k'}|D^N_{j'k'}M_b f(v_1,v_2)|^2dv\bigg)^{1/2} \\
     &\le C\|f\|_{L^2}\|g\|_{L^2}.
\endalign$$

To estimate $J_{12}$, let
$$K_{12}(z_1,u_1)=\int_{\Bbb R^m\times \Bbb R^m}
        \big[{}^t\overline D_{k}(y_2,z_2)-{}^t\overline D_{k}(y_2,v_2)\big]b_2(z_2)
   \widetilde K^1(z_1,u_1)(z_2,u_2)b_2(u_2)\overline D_{{k'}}(u_2,v_2) dz_2du_2$$
and $K_{12}$ be the operator associated with the kernel $K_{12}(z_1,u_1)$.
We claim that $K_{12}$ is a Calder\'on-Zygmund operator for $k'\ge k$ with
$$\|K_{12}\|_{CZ}\le C2^{-(k'-k)\varepsilon'}
     \frac{2^{-k\varepsilon'}}{(2^{-k}+|y_2-v_2|)^{m+\varepsilon'}}\qquad\text{for $0<\varepsilon'<\varepsilon$}.$$
We assume the claim for the moment and postpone its proof in Theorem 5. Rewrite $J_{12}$ as follows.
$$\allowdisplaybreaks\align
J_{12}
&=\sum_{j,k}\sum_{\Sb j'\ge j \\ k'\ge k \endSb}\int\int_{\Bbb R^n\times \Bbb R^n}
    \int M_bD^N_{jk}M_b g(y_1,y_2){}^t\!D_{j}(y_1,v_1)b_1(z_1)K_{12}(z_1,u_1)b_1(u_1)  \\
 &\hskip3cm \times   D_{{j'}}(u_1,v_1)M_bD^N_{j'k'}M_bf(v_1,v_2)    dydz_1du_1dv \\
 &=\sum_{j',k}\sum_{k'\ge k}\int\int_{\Bbb R^m}
              P_{j'}\overline D^N_{k}M_bg(v_1,y_2)b_2(y_2){}^t\!D_{{j'}}M_{b_1}({}^t\!K_{12}b_1)(v_1)M_b
             D^N_{j'k'}M_bf(v_1,v_2)   dvdy_2,
\endalign$$
where $P_{j'}(x_1,v_1):=\sum_{j\le j'}\int_{\Bbb R^n} D_{j}(x_1,y_1)b_1(y_1)D^N_{j}(y_1,v_1)dy_1$.
We use Schwarz's inequality twice and obtain
$$\align
|J_{12}|
&\le C\sum_{k}\sum_{k'\ge k}\int_{\Bbb R^m\times \Bbb R^m}
\bigg(\sum_{j'}\int_{\Bbb R^n} |P_{j'}\overline
D^N_{k}M_bg(v_1,y_2)|^2
          |{}^t\!D_{{j'}}M_{b_1}({}^t\!K_{12}b_1)(v_1)|^2 dv_1\bigg)^{1/2}\\
          &\hskip4cm \times   \bigg(\sum_{j'}\int_{\Bbb R^n} | D^N_{j'k'}M_bf(v_1,v_2)|^2 dv_1\bigg)^{1/2}dy_2dv_2. \\
\endalign$$
We first check that $P_{j'}(x_1,v_1)$ is bounded by a Poisson kernel.
Set $${\Cal P}_{j'}(x_1,v_1):=\sum_{-\infty<j\le j'}\int_{\Bbb R^n}
         \Cal V_N^{-1} D_{j}(\cdot,y_1)(x_1)b_1(y_1)D^N_{j}(y_1,v_1)dy_1$$ and
$${\Cal Q}_{j'}(x_1,v_1):=\sum_{j'< j<\infty}\int_{\Bbb R^n}
 \Cal V_N^{-1} D_{j}(\cdot,y_1)(x_1)b_1(y_1)D^N_{j}(y_1,v_1)dy_1,$$
where $\Cal V_N=\sum_{j\in \Bbb Z} D_jM_{b_1}D_j^NM_{b_1}$ is given in (1.2).
Then
$$f(x_1)=\sum_{j=-\infty}^\infty \Cal V_N^{-1}D_{j}M_{b_1}D^N_jM_{b_1}f(x_1)
={\Cal P}_{j'}M_{b_1}f(x_1)+{\Cal Q}_{j'}M_{b_1}f(x_1) \tag3.1$$
in $L^2(\Bbb R^n)$. For $|x_1-v_1|\le 2^{-j'}$,
$$|P_{j'}(x_1,v_1)|\le C\sum_{-\infty\le j\le j'} \frac{2^{jn}}{(1+2^j|x_1-v_1|)^{n+\varepsilon}}
     \le C\frac{2^{j'n}}{(1+2^{j'}|x_1-v_1|)^{n+\varepsilon}}.$$
Equality (3.1) shows ${\Cal P}_{j'}(x_1,v_1)b_1(v_1)+{\Cal Q}_{j'}(x_1,v_1)b_1(v_1)=\delta(x_1-v_1)$,
   the Dirac $\delta$-function. If $|x_1-v_1|>2^{-j'}$, then
$$|{\Cal P}_{j'}(x_1,v_1)|\le C|{\Cal Q}_{j'}(x_1,v_1)|
\le C\sum_{j'< j<\infty} \frac{2^{jn}}{(1+2^j|x_1-v_1|)^{n+\varepsilon}}
     \le C\frac{2^{j'n}}{(1+2^{j'}|x_1-v_1|)^{n+\varepsilon}}.$$
By [H, Theorem 2.8],
$$\align |P_{j'}(x_1,v_1)|
&=|\Cal V_N {\Cal P}_{j'}(\cdot,v_1)(x_1)|\le C\frac{2^{j'n}}{(1+2^{j'}|x_1-v_1|)^{n+\varepsilon}}.
\endalign$$
Therefore,
$$| P_{j'}(x_1,v_1)|\le C\frac{2^{j'n}}{(1+2^{j'}|x_1-v_1|)^{n+\varepsilon}}
     \qquad\text{for all $x_1,v_1\in \Bbb R^n$}.\tag 3.2$$

For $k'\ge k$, we assume $K_{12}$ to be bounded on $L^2$ with $\|K_{12}\|_{CZ}\le
C2^{-(k'-k)\varepsilon'}
         \frac{2^{-k\varepsilon'}}{(2^{-k}+|y_2-v_2|)^{m+\varepsilon'}}$. The $Tb$ theorem shows
$\!K_{12}b_1\in BMO(\Bbb R^n)$. By (3.2) and Carleson measure condition on $\Bbb R^n$,
$$\align
|J_{12}|
      &\le  C\sum_{k}\sum_{k'\ge k}\int_{\Bbb R^m\times \Bbb R^m}
       \bigg(\int_{\Bbb R^n} |\overline D^N_{k}M_{b_2}g(v_1,y_2)|^2dv_1\bigg)^{1/2}2^{-(k'-k)\varepsilon'}
         \frac{2^{-k\varepsilon'}}{(2^{-k}+|y_2-v_2|)^{m+\varepsilon'}} \\
      &\hskip5cm \times  \bigg(\int_{\Bbb R^n} |\overline D^N_{k'}M_{b_2}f(v_1,v_2)|^2
      dv_1\bigg)^{1/2}dy_2dv_2.
\endalign$$
Using Schwarz's inequality and $L^2$-boundedness of $\Cal G$-function, we have
$$\allowdisplaybreaks\align
|J_{12}| & \le C \bigg(\sum_{k}\sum_{k'\ge k}\int_{\Bbb R^m}\int
|\overline D^N_{k}M_{b_2}g(v_1,y_2)|^22^{-(k'-k)\varepsilon'}
         \frac{2^{-k\varepsilon'}}{(2^{-k}+|y_2-v_2|)^{m+\varepsilon'}}dy_2dv \bigg)^{1/2}\\
   &\qquad\times        \bigg(\sum_{k}\sum_{k'\ge k}\int_{\Bbb R^m}\int  |\overline D^N_{k'}M_{b_2}f(v_1,v_2)|^2 2^{-(k'-k)\varepsilon'}
         \frac{2^{-k\varepsilon'}}{(2^{-k}+|y_2-v_2|)^{m+\varepsilon'}}dy_2dv \bigg)^{1/2}\\
   &\le C\bigg( \int_{\Bbb R^n\times \Bbb R^m}
       \sum_{k}|\overline D^N_{k}M_{b_2}g(v_1,y_2)|^2dv_1dy_2\bigg)^{1/2}
       \bigg(\int \sum_{k'}|\overline D^N_{k'}M_{b_2}f(v_1,v_2)|^2 dv \bigg)^{1/2}\\
   &\le C\|f\|_{L^2}\|g\|_{L^2}.
 \endalign$$
The estimate of $J_{13}$ is similar to $J_{12}$, and we leave details to
readers. For $J_{14}$, we write
$$\allowdisplaybreaks\align
J_{14}
&=\sum_{j,k}\sum_{\Sb j'\ge j \\ k'\ge k \endSb}
          \iiint M_bD^N_{jk}M_b g(y_1,y_2){}^t\!D_j(y_1,v_1){}^t\overline D_k(y_2,v_2)\\
&\qquad\qquad\times {}^tTb(u_1,u_2)b(u_1,u_2)D_{j'}(u_1,v_1)\overline D_{k'}(u_2,v_2)
           M_bD^N_{j'k'}M_bf(v_1,v_2)dudydv\\
 &=\sum_{j',k'} \iint {}^tTb(u_1,u_2) b_1(u_1)D_{j'}(u_1,v_1)b_2(u_2)\overline D_{k'}(u_2,v_2)\\
    &\hskip4cm\times    M_bD^N_{j'k'}M_bf(v_1,v_2) P_{j'}\overline P_{k'}M_bg(v_1,v_2)dudv\\
    &=\sum_{j',k'}  \int {}^t\!D_{j'}{}^t\overline D_{k'}M_b({}^tTb)(v_1,v_2)M_bD^N_{j'k'}M_bf(v_1,v_2)
                P_{j'}\overline P_{k'}M_bg(v_1,v_2)dv,
\endalign$$
where $\overline P_{k'}(x_2,v_2):=\sum\limits_{k\le k'}\int_{\Bbb
R^m} \overline D_{k}(x_2,y_2)$ $b_2(y_2)\overline
D^N_{k}(y_2,v_2)dy_2$. By ${}^tTb\in BMO$ and  Carleson measure condition on $\Bbb R^n\times \Bbb R^m$,
$$\allowdisplaybreaks\align
|J_{14}|&\le C\Big(\int \sum_{j',k'}
      |{}^t\!D_{j'}{}^t\overline D_{k'}M_b({}^tTb)(v_1,v_2)|^2| P_{j'}\overline
      P_{k'}M_bg(v_1,v_2)|^2dv\Big)^{1/2}\\
&\hskip3cm\times\Big(\int \sum_{j',k'}
      |D^N_{j'k'}M_bf(v_1,v_2)|^2dv\Big)^{1/2}\\
&\le C\|{}^tTb\|_{BMO}\|f\|_{L^2}\|g\|_{L^2}\le C\|f\|_{L^2}\|g\|_{L^2}.
\endalign$$

To estimate $J_2$, we also need to use the almost
orthogonality argument, and write
$$\allowdisplaybreaks\align
J_2
&=\sum_{j,k}\sum_{\Sb j'\ge j \\ k'< k \endSb} \bigg\langle
M_bD^N_{jk}M_b g, \iiint
 \big[{}^t\!D_j(\cdot,z_1)-{}^t\!D_j(\cdot,v_1)\big]b_1(z_1)
        {}^t\overline D_k(\cdot,z_2)b_2(z_2)K(z_1,z_2,u_1,u_2)  \\
&\qquad\qquad\qquad \times b_1(u_1)D_{{j'}}(u_1,v_1)
             b_2(u_2)\big[\overline D_{{k'}}(u_2,v_2)-\overline
             D_{{k'}}(\cdot,v_2)\big]M_bD^N_{j'k'}M_bf(v_1,v_2)dzdudv \bigg\rangle \\
&\ + \sum_{j,k}\sum_{\Sb j'\ge j \\ k'< k \endSb} \bigg\langle
M_bD^N_{jk}M_b g, \iiint  {}^t\!D_j(\cdot,v_1)b_1(z_1)
        {}^t\overline D_k(\cdot,z_2)b_2(z_2) K(z_1,z_2,u_1,u_2) \\
&\qquad\qquad\qquad \times b_1(u_1)D_{{j'}}(u_1,v_1)
             b_2(u_2)\big[\overline D_{{k'}}(u_2,v_2)-\overline
             D_{{k'}}(\cdot,v_2)\big]M_bD^N_{j'k'}M_bf(v_1,v_2)dzdudv \bigg\rangle\\
&\ +\sum_{j,k}\sum_{\Sb j'\ge j \\ k'< k \endSb} \bigg\langle
            M_bD^N_{jk}M_b g, \iiint \big[{}^t\!D_j(\cdot,z_1)-{}^t\!D_j(\cdot,v_1)\big]b_1(z_1)
        {}^t\overline D_k(\cdot,z_2)b_2(z_2) K(z_1,z_2,u_1,u_2) \\
&\qquad\qquad\qquad \times b_1(u_1)D_{{j'}}(u_1,v_1)
             b_2(u_2)\overline
             D_{{k'}}(\cdot,v_2)M_bD^N_{j'k'}M_bf(v_1,v_2)dzdudv \bigg\rangle\\
&\ +\sum_{j,k}\sum_{\Sb j'\ge j \\ k'< k \endSb} \bigg\langle
M_bD^N_{jk}M_b g, \iiint {}^t\!D_j(\cdot,v_1)b_1(z_1)
        {}^t\overline D_k(\cdot,z_2)b_2(z_2) K(z_1,z_2,u_1,u_2) \\
&\qquad\qquad\qquad \times b_1(u_1)D_{{j'}}(u_1,v_1)
             b_2(u_2)\overline D_{{k'}}(\cdot,v_2)M_bD^N_{j'k'}M_bf(v_1,v_2)dzdudv \bigg\rangle\\
&:=J_{21}+J_{22}+J_{23}+J_{24}
\endalign$$
By using the almost orthogonality, a similar argument to the estimate of $J_{11}$ shows
$|J_{21}|\le C\|f\|_{L^2}\|g\|_{L^2}$. To estimate $J_{22}$, we set
$$K_{22}(z_1,u_1):=\int_{\Bbb R^m\times\Bbb R^m}  {}^t\overline D_{k}(y_2,z_2)b_2(z_2)
                                \widetilde K^1(z_1,u_1)(z_2,u_2)b_2(u_2)
                           \big[\overline D_{{k'}}(u_2,v_2)-\overline D_{{k'}}(y_2,v_2)\big]dz_2du_2.$$
By Theorem 5 below, the operator $K_{22}$ associated with the kernel
            $K_{22}(z_1,u_1)$ is a Calder\'on-Zygmund operator for $k'\le k$ and
$$\|K_{22}\|_{CZ}\le C2^{-(k-k')\varepsilon'}
         \frac{2^{-k'\varepsilon'}}{(2^{-k'}+|y_2-v_2|)^{m+\varepsilon'}}
                                  \qquad\text{for $0<\varepsilon'<\varepsilon$}.$$
Thus,
$$\align
J_{22}
&=\sum_{j,k'}\sum_{\Sb j'\ge j \\ k'< k \endSb}
     \iint_{\Bbb R^n\times \Bbb R^n}\int M_bD^N_{jk}M_b g(y_1,y_2)
         {}^t\!D_{j}(y_1,v_1)b_1(z_1)K_{22}(z_1,u_1) \\
&\hskip3.5cm \times b_1(u_1)D_{{j'}}(u_1,v_1)M_bD^N_{{j'k'}}M_bf(v_1,v_2)dydz_1du_1dv.
\endalign$$
The estimate of $J_{22}$, $|J_{22}|\le C\|f\|_{L^2}\|g\|_{L^2}$,
             is done by using the same argument as the proof of $J_{12}$.
Since $J_{23}$ is a symmetric case of $J_{22}$,
             it remains $J_{24}$ for considering. To do this, we write
$$\allowdisplaybreaks\align
J_{24}
&=\sum_{j,k}\sum_{\Sb j'\ge j \\ k'< k \endSb}\int M_bD^N_{jk}M_b
        g(y_1,y_2) \int_{\Bbb R^n\times \Bbb R^m}\int {}^t\!D_j(y_1,v_1)
        {}^t\overline D_k(y_2,z_2)b_2(z_2) \widetilde Tb(u_1,z_2) \\
&\qquad\qquad\qquad \times b_1(u_1)D_{{j'}}(u_1,v_1)
             \overline D_{{k'}}(y_2,v_2)M_bD^N_{j'k'}M_bf(v_1,v_2)dydu_1dz_2dv\\
&=\sum_{j,k}\sum_{\Sb j'\ge j \\ k'< k \endSb} \Big\langle bg,
             \iiint_{\Bbb R^n\times \Bbb R^m}
     {}^t\!D^N_j(\cdot,y_1)b_1(y_1){}^t\!D_j(y_1,v_1)
     {}^t\overline D^N_k(\cdot,y_2)b_2(y_2){}^t\overline D_k(y_2,z_2)b_2(z_2)  \\
&\qquad\qquad\qquad \times \widetilde Tb(u_1,z_2)b_1(u_1)D_{{j'}}(u_1,v_1)
             \overline D_{{k'}}(y_2,v_2)M_bD^N_{j'k'}M_bf(v_1,v_2)dydvdu_1dz_2\Big\rangle \\
&=\Big\langle  bg, \sum_{j',k}\int_{\Bbb R^m\times\Bbb R^m}\int_{\Bbb R^n\times \Bbb R^n}\int
             {}^t\!P_{j'}(\cdot,v_1){}^t\overline D^N_{k}(\cdot,y_2)b_2(y_2)
                {}^t\overline D_{k}(y_2,z_2)b_2(z_2) \widetilde Tb(u_1,z_2) b_1(u_1) \\
 &\qquad\qquad\quad \times      D_{{j'}}(u_1,v_1)b_1(v_1)D^N_{{j'}}(v_1,w_1)
             \overline P_{k+1}(y_2,w_2)(bf)(w_1,w_2) dy_2dz_2du_1dv_1dw\Big\rangle\\
&:=\big\langle bg, S_{\widetilde Tb}(f)  \big\rangle.
\endalign$$
If we have, for $h\in L^\infty(\Bbb R^n\times \Bbb R^m)$,
$$\|S_{\widetilde Tb}(h)\|_{BMO}\le C\| \widetilde Tb\|_{BMO}\|h\|_{L^\infty}\tag 3.3$$
and
 $$\|{}^t\!S_{\widetilde Tb}(h)\|_{BMO}\le  C\| \widetilde Tb\|_{BMO}\|h\|_{L^\infty},\tag 3.4$$
then, for $h\in H^1(\Bbb R^n\times \Bbb R^m)$,
$$\|S_{\widetilde Tb}(h)\|_{L^1}=\big|\big\langle S_{\widetilde Tb}(h),
                        \operatorname{sgn}\{S_{\widetilde Tb}(h)\}\big\rangle\big|
 = \big|\big\langle h, {}^t\!S_{\widetilde Tb}(\operatorname{sgn}\{S_{\widetilde Tb}(h)\})\big\rangle\big|
 \le C\| \widetilde Tb\|_{BMO}\|h\|_{H^1}.\tag3.5$$
By the assumption $\widetilde Tb\in BMO$ and using interpolation between (3.3) and (3.5), we obtain
$$|J_{24}|=\big|\big\langle bg,  S_{\widetilde Tb}(f)\big\rangle\big|
            \le C\|S_{\widetilde Tb}(f)\|_{L^2}\|g\|_{L^2}\le C\|f\|_{L^2}\|g\|_{L^2}.$$
Since the proofs of (3.3) and (3.4) are similar, we check (3.3) only.
For $h\in L^\infty(\Bbb R^n\times \Bbb R^m)$, consider the operator $\Bbb S_h$
on $BMO(\Bbb R^n\times \Bbb R^m)$ given by $\Bbb S_h(\phi)=S_\phi(h)$.
We will show that $\Bbb S_h$ is bounded on $BMO(\Bbb R^n\times \Bbb R^m)$.
Observe the transport ${}^t\Bbb S_h$ which is
$$\allowdisplaybreaks\align
{}^t\Bbb S_h(\phi)(u_1,z_2)
 &=\sum_{j',k} \int_{\Bbb  R^n\times \Bbb R^m}
 b_1(u_1)D_{{j'}}(u_1,v_1)b_1(v_1)b_2(y_2){}^t\!\overline D_{k}(y_2,z_2)b_2(z_2) \\
 &\qquad\qquad\qquad \times       P_{j'}\overline D^N_{k}\phi(v_1,y_2)D^N_{{j'}}
                      \overline P_{k+1}M_b h(v_1,y_2)dv_1dy_2.\\
\endalign$$
First, we show that
 $\|{}^t\Bbb S_h(\phi)\|_{L^2}\le C\|\phi\|_{L^2}\|h\|_{L^\infty}$
provided $h\in L^\infty$. It follows from $\overline P_{k+1} bh(\cdot,y_2)\in L^\infty\subset BMO$
and the Carleson measure condition on $\Bbb R^n$ that
$$\allowdisplaybreaks\align
&\|{}^t\Bbb S_h(\phi)\|_{L^2}\\
  &=\sup_{\|q\|_{L^2}\le 1} \sum_{j',k}\int_{\Bbb R^n\times\Bbb R^m}
  P_{j'}\overline D^N_k\phi(v_1,y_2)b_2(y_2)
          {}^t\!D_{j'}{}^t\overline D_kM_bq(v_1,y_2) b_1(v_1)D^N_{j'}\overline P_{k+1}M_b h(v_1,y_2)dv_1dy_2 \\
 &\le C\sup_{\|q\|_{L^2}\le 1} \|q\|_{L^2} \bigg\{\int_{\Bbb R^n\times\Bbb R^m}
            \sum_{j',k} |P_{j'}\overline D^N_k\phi(v_1,y_2)|^2
                      |D^N_{j'}\overline P_{k+1}M_b h(v_1,y_2)|^2dv_1dy_2\bigg\}^{1/2} \\
 &\le C\|h\|_{L^\infty} \bigg\{\int_{\Bbb R^n\times\Bbb R^m}
           \sum_k |\overline D^N_k\phi(v_1,y_2)|^2 dv_1dy_2\bigg\}^{1/2}\\
 &\le C\|h\|_{L^\infty}\| \phi\|_{L^2},
\endalign$$
where the last inequality is obtained by the Littlewood-Paley estimate on $\Bbb R^m$. Using the Carleson measure condition, property (3.2), and
$$\int_{\Bbb R^n} \frac{2^{-j\varepsilon}}{(2^{-j}+|x-v|)^{n+\varepsilon}}
            \frac{2^{-j\varepsilon}}{(2^{-j}+|v-y|)^{n+\varepsilon}}dv
        \le C\frac{2^{-j\varepsilon}}{(2^{-j}+|x-y|)^{n+\varepsilon}},$$
we may obtain that the kernel $\Bbb K(x_1,x_2,y_1,y_2):=b^{-1}(x_1,x_2){}^t\Bbb S_h(x_1,x_2,y_1,y_2)$ of
$M_{b^{-1}}{}^t\Bbb S_h$ satisfies conditions $(A_1)$ and $(A_2)$ since
$$\align
{}^t\Bbb S_h(x_1,x_2,y_1,y_2) &=\sum_{j',k} \int
b_1(x_1)D_{{j'}}(x_1,v_1)b_2(x_2)\overline
D_{k}(x_2,v_2)b_1(v_1)b_2(v_2)\\
 &\qquad\qquad \times    P_{j'}(v_1,y_1)\overline D^N_{k}(v_2,y_2)
   D^N_{{j'}}\overline P_{k+1}M_bh(v_1,v_2)dv.
\endalign$$
Let $$C^\infty_{0,0}(\Bbb R^n) :=\bigg\{ \psi \in C^\infty(\Bbb R^n) : \psi\
\text{has a compact support}\ \text{and}\ \int_{\Bbb R^n} \psi=0 \bigg\}.$$
The properties of $D_{j'}(x_1,v_1)$ and $\overline D_k(x_2,v_2)$ give
$$\cases\displaystyle \int_{\Bbb R^n}\int
         b(x_1,x_2)\Bbb K(x_1,x_2,y_1,y_2)\varphi^1(y_1)\varphi^2(y_2)dx_1dy=0,\qquad x_2\in \Bbb R^m,\\
     \displaystyle \int_{\Bbb R^m}\int b(x_1,x_2)\Bbb K(x_1,x_2,y_1,y_2)\varphi^1(y_1)\varphi^2(y_2)dx_2dy=0,
         \qquad x_1\in \Bbb R^n,
\endcases$$
for all $\varphi^1\in C^\infty_{0,0}(\Bbb R^n)$ and $\varphi^2\in C^\infty_{0,0}(\Bbb R^m)$.
The same proof of [HLLL, Theorem 1] implies
$$\|{}^t{\Bbb S_h}(\phi)\|_{H^1}\le
       C\|\phi\|_{H^1}\|h\|_{L^\infty}.$$
By duality, $\|S_{\widetilde Tb}(h)\|_{BMO}
   =\|\Bbb S_h(\widetilde Tb)\|_{BMO}\le C\| \widetilde Tb\|_{BMO}\|h\|_{L^\infty}$
and (3.3) follows.

\vskip0.2cm
To complete the proof of Theorem 1, we still need to show the following

\proclaim{Theorem 5} Let $b$ be a para-accretive function defined on
$\Bbb R^m$ and $\{\overline S_k\}_{k\in \Bbb Z}$
 be an approximation to the identity associated to $b$ with
regularity exponent $\varepsilon.$ Set $\overline D_k=\overline S_k-\overline S_{k-1}$.
 Let $T$ be a generalized singular integral operator associated to a kernel $K(z_1,z_2,u_1,u_2)$ with
regularity exponent $\varepsilon.$ For $k'\ge k$, define
$$\Cal K(z_1,u_1)=\int_{\Bbb R^m\times\Bbb R^m} \big[\overline D_k(y_2,z_2)-\overline D_k(y_2,v_2)\big]
                                             b(z_2)K(z_1,z_2,u_1,u_2)b(u_2)\overline D_{k'}(u_2,v_2)dz_2du_2.$$
Then the operator $\Cal K$ given by
$$\langle \Cal Kf, g \rangle =\iint  g(x)\Cal K(x,y)f(y)dxdy,
                          \qquad\operatorname{supp}(f)\cap \operatorname{supp}(g)=\varnothing,$$
is a Calder\'on-Zygmund operator and satisfies, for $0<\varepsilon'<\varepsilon$,
$$\|\Cal K\|_{CZ} \le C2^{-(k'-k)\varepsilon'}
                  \frac{2^{-k\varepsilon'}}{(2^{-k}+|y_2-v_2|)^{m+\varepsilon'}}.$$
\endproclaim

\demo{Proof}
First we show that $\Cal K(z_1,u_1)$ is a Calder\'on-Zygmund kernel.
For fixed $z_1$ and $u_1$,\break ${\widetilde K}^1(z_1,u_1)(z_2,u_2)=K(z_1,z_2,u_1,u_2)$ is a Calder\'on-Zygmund
kernel on $\Bbb R^m\times \Bbb R^m\backslash \{z_2=u_2\}$ with the norm
$\|{\widetilde K}^1(z_1,u_1)\|_{CZ}\le C|z_1-u_1|^{-n}.$
By the almost orthogonality estimate for the kernel ${\widetilde K}^1(z_1,u_1)(z_2,u_2)$,
$$\aligned |\Cal K(z_1,u_1)|
&\le C2^{-(k'-k)\varepsilon'}
                 \frac{2^{-k\varepsilon'}}{(2^{-k}+|y_2-v_2|)^{m+\varepsilon'}}\big\|{\widetilde K}^1(z_1,u_1)\big\|_{CZ}\\
&\le  C2^{-(k'-k)\varepsilon'}\frac{2^{-k\varepsilon'}}{(2^{-k}+|y_2-v_2|)^{m+\varepsilon'}}|z_1-u_1|^{-n}.
 \endaligned\tag 3.6 $$
It is easy to see
$$\align
\Cal K(z_1,u_1)-\Cal K(z_1,u'_1)&=\int_{\Bbb R^m\times \Bbb R^m}\hskip-0.2cm
                           \big[\overline D_k(y_2,z_2)-\overline D_k(y_2,v_2)\big]b(z_2)\\
&\qquad\times \big[K(z_1,z_2,u_1,u_2)-K(z_1,z_2,u'_1,u_2)\big]b(u_2)\overline D_{k'}(u_2,v_2)dz_2du_2.
\endalign$$
Note that ${\widetilde K}^1(z_1,u_1)-{\widetilde K}^1(z_1,u'_1)$
              is a Calder\'on-Zygmund operator on $\Bbb R^m$ with
$\big\|{\widetilde K}^1(z_1,u_1)-{\widetilde K}^1(z_1,u'_1)\big\|_{CZ}\le C |u_1-u'_1|^{\varepsilon}
      |z_1-u_1|^{-n-\varepsilon}$ for $|u_1-u'_1|\le |z_1-u_1|/2.$
The same argument as (3.6) gives
$$|\Cal K(z_1,u_1)-\Cal K(z_1,u'_1)|\le C2^{-(k'-k)\varepsilon'}
       \frac{2^{-k\varepsilon'}}{(2^{-k}+|y_2-v_2|)^{m+\varepsilon'}}
      |u_1-u'_1|^{\varepsilon}
      |z_1-u_1|^{-n-\varepsilon}$$
for $|u_1-u'_1|\le |z_1-u_1|/2$. Similarly,
$$|\Cal K(z_1,u_1)-\Cal K(z'_1,u_1)|
  \le C 2^{-(k'-k)\varepsilon'}
                   \frac{2^{-k\varepsilon'}}{(2^{-k}+|y_2-v_2|)^{m+\varepsilon'}}
         |z_1-z'_1|^{\varepsilon}|z_1-u_1|^{-n-\varepsilon}$$
for $|z_1-z'_1|\le|z_1-u_1|/2$.
Hence, $\Cal K(z_1,u_1)$ is a Calder\'on-Zygmund kernel with
$$|\Cal K|_{CZ} \le C2^{-(k'-k)\varepsilon'}
                    \frac{2^{-k\varepsilon'}}{(2^{-k}+|y_2-v_2|)^{m+\varepsilon'}}.$$

Next we show that $\Cal K$ is bounded on $L^2(\Bbb R^n)$ with
$\|\Cal Kf\|_{L^2}\le C2^{-(k'-k)\varepsilon'}
           \frac{2^{-k\varepsilon'}}{(2^{-k}+|y_2-v_2|)^{m+\varepsilon'}}\|f\|_{L^2}$.
By duality,
$$\align
\|\Cal Kf\|_{L^2}
&=\sup_{\|h\|_{L^2}\le1} \int_{\Bbb R^m\times \Bbb R^m}\int_{\Bbb R^n\times \Bbb R^n} h(z_1)
                   \big[\overline D_k(y_2,z_2)-\overline D_k(y_2,v_2)\big] \\
&\hskip4cm\times            b(z_2)K(z_1,z_2,u_1,u_2)b(u_2)\overline D_{k'}(u_2,v_2)
              f(u_1)dz_2du_2dz_1du_1\\
&=\sup_{\|h\|_{L^2}\le1} \bigg(\int_{|z_2-v_2|> 8c2^{-k}}+\int_{|z_2-v_2|\le 8c2^{-k}}\bigg)
               \int_{\Bbb R^m}\int_{\Bbb R^n\times \Bbb R^n} h(z_1)
                   \big[\overline D_k(y_2,z_2)-\overline D_k(y_2,v_2)\big] \\
&\hskip4cm\times            b(z_2)K(z_1,z_2,u_1,u_2)b(u_2)\overline D_{k'}(u_2,v_2)
              f(u_1)dz_2du_2dz_1du_1\\
&:= I_1+I_2,
\endalign$$
where the constant $c$ satisfies item (ii) in the definition od $\{\overline S_k\}$.
For $I_1$, we use the cancellation property of $\overline D_{k'}$ to get
$$\align
I_1
&=\sup_{\|h\|_{L^2}\le1} \int_{|z_2-v_2|> 8c2^{-k}}\int_{\Bbb R^m}\int_{\Bbb R^n\times \Bbb R^n}
               h(z_1)\big[\overline D_k(y_2,z_2)-\overline D_k(y_2,v_2)\big] \\
&\qquad\qquad\times                  b(z_2)\big[K(z_1,z_2,u_1,u_2)-K(z_1,z_2,u_1,v_2)\big]
           b(u_2)\overline D_{k'}(u_2,v_2)f(u_1)dz_1du_1dz_2du_2\\
&=\sup_{\|h\|_{L^2}\le1}  \int_{|z_2-v_2|> 8c2^{-k}}\int_{\Bbb R^m} \bigg(\int_{\Bbb R^n\times \Bbb R^n}
           \!\!\! h(z_1)\big[{\widetilde K}^2(z_2,u_2)-{\widetilde K}^2(z_2,v_2)\big](z_1,u_1)f(u_1)dz_1du_1\bigg) \\
     &\hskip5cm \times        \big[\overline D_k(y_2,z_2)-\overline D_k(y_2,v_2)\big]
                                b(z_2)b(u_2)\overline D_{k'}(u_2,v_2)dz_2du_2.\\
\endalign$$
In the above integral, if $\overline D_{k'}(u_2,v_2)\ne0$, then $|u_2-v_2|\le 2c2^{-k'}$.
Since $k'\ge k$ and $|z_2-u_2|\ge |z_2-v_2|-|u_2-v_2|\ge 6c2^{-k}$, we have $|u_2-v_2|\le |z_2-u_2|/2$. Hence,
$$\align
I_1
&\le C\|f\|_{L^2}\int_{\Bbb R^m\times \Bbb R^m}\frac {|u_2-v_2|^\varepsilon}{|z_2-u_2|^{m+\varepsilon}}
              \big[|\overline D_k(y_2,z_2)|+|\overline D_k(y_2,v_2)|\big]|\overline D_{k'}(u_2,v_2)|dz_2du_2\\
&\le C2^{-(k'-k)\varepsilon}\|f\|_{L^2} \int_{\Bbb R^m\times \Bbb R^m}
   \frac {2^{-k\varepsilon}}{(2^{-k\varepsilon}+|z_2-u_2|)^{m+\varepsilon}}|\overline D_k(y_2,z_2)||\overline D_{k'}(u_2,v_2)|dz_2du_2\\
&\qquad+ C2^{-(k'-k)\varepsilon}\|f\|_{L^2} |\overline D_k(y_2,v_2)|\int_{\Bbb R^m\times \Bbb R^m}
   \frac {2^{-k\varepsilon}}{(2^{-k\varepsilon}+|z_2-u_2|)^{m+\varepsilon}}|\overline D_{k'}(u_2,v_2)|dz_2du_2\\
&\le C2^{-(k'-k)\varepsilon'}\frac{2^{-k\varepsilon'}}{(2^{-k}+|y_2-v_2|)^{m+\varepsilon'}}\|f\|_{L^2}.
\endalign$$

For $I_2$, by the condition on the support of $\overline D_k$, we have $|y_2-v_2|\le 10C2^{-k}$.
Let $\eta_0\in C^\infty(\Bbb R^m)$ be 1 on the unit ball and 0
outside the ball $B(0,2)$ and set $\eta_1 = 1- \eta_0$.
$$\align
I_2&=\sup_{\|h\|_{L^2}\le1} \int_{|z_2-v_2|\le 8c2^{-k}}
               \int_{\Bbb R^m}\int_{\Bbb R^n\times \Bbb R^n} h(z_1)
                   \big[\overline D_k(y_2,z_2)-\overline D_k(y_2,v_2)\big] \\
&\hskip4cm\times            b(z_2)K(z_1,z_2,u_1,u_2)b(u_2)\overline D_{k'}(u_2,v_2)
              f(u_1)dz_2du_2dz_1du_1\\
&=\sup_{\|h\|_{L^2}\le1} \int_{|z_2-v_2|\le 8c2^{-k}}
               \int_{\Bbb R^m}\int_{\Bbb R^n\times \Bbb R^n} h(z_1)
                    \eta_0\Big(\frac{z_2-v_2}{c2^{-k'+2}}\Big)\big[\overline D_k(y_2,z_2)-\overline D_k(y_2,v_2)\big] \\
&\hskip4cm\times            b(z_2)K(z_1,z_2,u_1,u_2)b(u_2)\overline D_{k'}(u_2,v_2)
              f(u_1)dz_2du_2dz_1du_1\\
&\qquad+ \sup_{\|h\|_{L^2}\le1} \int_{|z_2-v_2|\le 8c2^{-k}}
               \int_{\Bbb R^m}\int_{\Bbb R^n\times \Bbb R^n} h(z_1)
                    \eta_1\Big(\frac{z_2-v_2}{c2^{-k'+2}}\Big)\big[\overline D_k(y_2,z_2)-\overline D_k(y_2,v_2)\big] \\
&\hskip4cm\times            b(z_2)K(z_1,z_2,u_1,u_2)b(u_2)\overline D_{k'}(u_2,v_2)
              f(u_1)dz_2du_2dz_1du_1\\
&:=I_{21}+I_{22}.
\endalign$$
To estimate $I_{21}$, we define $\Cal K^0(z_1,u_1)$ by
$$\align
\Cal K^0(z_1,u_1)&=\int_{|z_2-v_2|\le 8c2^{-k}}
               \int_{\Bbb R^m} \eta_0\Big(\frac{z_2-v_2}{c2^{-k'+2}}\Big)
         \big[\overline D_k(y_2,z_2)-\overline D_k(y_2,v_2)\big]\\
&\qquad\qquad\times b(z_2)
           K(z_1,z_2,u_1,u_2)b(u_2)\overline D_{k'}(u_2,v_2)dz_2du_2.
\endalign$$
Let $F(z_2):=\eta_0\big(\frac{z_2-v_2}{c2^{-k'+2}}\big)
   (\overline D_k(y_2,z_2)-\overline D_k(y_2,v_2))$ and $G(u_2):=\overline D_{k'}(u_2,v_2)$.
Weak boundedness property shows that
$$\aligned I_{21}
&=\sup_{\|h\|_{L^2}\le1}\int_{\Bbb R^n\times\Bbb R^n} h(z_1)\Cal K^0(z_1,u_1)f(u_1)dz_1du_1  \\
&=\sup_{\|h\|_{L^2}\le1}\int_{\Bbb R^n\times\Bbb R^n}\int_{|z_2-v_2|\le 8c2^{-k}}
               \int_{\Bbb R^m} h(z_1)\eta_0\Big(\frac{z_2-v_2}{c2^{-k'+2}}\Big)
         \big[\overline D_k(y_2,z_2)-\overline D_k(y_2,v_2)\big]\\
&\qquad\times b(z_2)
           K(z_1,z_2,u_1,u_2)b(u_2)\overline D_{k'}(u_2,v_2)dz_2du_2f(u_1)dz_1du_1\\
&\le \sup_{\|h\|_{L^2}\le1} C2^{-k'(m+2\eta)}\|F\|_\eta\|G\|_\eta\|h\|_{L^2}\|f\|_{L^2}\\
&\le C2^{-(k'-k)\varepsilon}2^{km}\|f\|_{L^2}  \\
&\le C2^{-(k'-k)\varepsilon}\frac{2^{-k\varepsilon'}}{(2^{-k}+|y_2-v_2|)^{m+\varepsilon'}}\|f\|_{L^2}.
\endaligned\tag 3.7$$
To estimate $I_{22}$, we use the cancellation property of $\overline D_{k'}$ and write
$$\align
I_{22}
&=\sup_{\|h\|_{L^2}\le1} \int_{|z_2-v_2|\le 8c2^{-k}}
               \int_{\Bbb R^m}\int_{\Bbb R^n\times \Bbb R^n} h(z_1)
                    \eta_1\Big(\frac{z_2-v_2}{c2^{-k'+2}}\Big)\big[\overline D_k(y_2,z_2)-\overline D_k(y_2,v_2)\big] \\
&\qquad\qquad\times            b(z_2)\big[K(z_1,z_2,u_1,u_2)-K(z_1,z_2,u_1,v_2)\big]b(u_2)\overline D_{k'}(u_2,v_2)
              f(u_1)dz_2du_2dz_1du_1
\endalign$$
By the conditions on the supports of $\eta_1$ and $\overline D_{k'}$, we have $|z_2-v_2|\ge 4c2^{-k'}$ and $|u_2-v_2|\le 2c2^{-k'}$.
This gives $|u_2-v_2|\le |z_2-v_2|/2$. Applying $(A_2)$, we obtain
$$\align
II&\le C\|f\|_{L^2} \int_{4c2^{-k'}<|z_2-v_2|\le 8c2^{-k}}\int_{|u_2-v_2|\le 2c2^{-k'}}
              |\overline D_k(y_2,z_2)-\overline D_k(y_2,v_2)|\\
&\hskip3cm \times \|\widetilde K^2(z_2,u_2)-\widetilde K^2(z_2,v_2)\|_{CZ}|\overline D_{k'}(u_2,v_2)|dz_2du_2 \\
&\le C\|f\|_{L^2}  \int_{4c2^{-k'}<|z_2-v_2|\le 8c2^{-k}}\int_{|u_2-v_2|\le 2c2^{-k'}}
 |z_2-v_2|^\varepsilon 2^{k(m+\varepsilon)}
            \frac{|u_2-v_2|^\varepsilon}{|z_2-v_2|^{m+\varepsilon}}2^{k'm}dz_2du_2 \\
&\le C\|f\|_{L^2} 2^{-(k'-k)\varepsilon}2^{km} \int_{4c2^{-k'}<|z_2-v_2|\le 8c2^{-k}}  \frac 1{|z_2-v_2|^m} dz_2\\
&\le C\|f\|_{L^2} 2^{-(k'-k)\varepsilon'}2^{km}.
\endalign$$
Since $|y_2-v_2|\le 10c2^{-k}$, we get
$II\le C2^{-(k'-k)\varepsilon'}\frac{2^{-k\varepsilon'}}{(2^{-k}+|y_2-v_2|)^{m+\varepsilon'}}\|f\|_{L^2}.$
Combining with (3.7), we finish the proof of Theorem 5. \hfill $\qed$
\enddemo

\vskip 0.3cm
\noindent
{\bf Note added in revision.}
Although the $T1$ and $Tb$ theorems were obtained earlier than two decades ago, it is highly desirable to acquire further insight the theory on product spaces.
After the authors submitted the manuscript to arXiv, they learned from Michael Lacey that several other authors did some extensions in this subject recently.
Pott and Villaroya [PV] presented a nice extension of Journe's $T1$ theorem, and
Martikainen [Ma] got a modern approach to things.
Ou [O] also obtained a $Tb$ theorem which includes the advances of the previous two papers.
In a different direction, Hyt\"onen and Martikainen [HM] showed a non-homogeneous $T1$ theorem for two parameters.
The authors are grateful to Michael Lacey for bringing these helpful references to their attention.

The research was initiated when the first author visited Taiwan in April 2009.
He acknowledges a financial support received from NCU Center for
Mathematics and Theoretic Physics and hospitality offered by the Department of Mathematics,
National Central University, Taiwan, Republic of China.

\noindent


\vskip 0.3in \widestnumber\key{HLLL} \Refs\nofrills{References}

\ref\key DJ
\by G. David and J.-L. Journ\'e
\paper A boundedness criterion for generalized Calder\'on-Zygmund operators
\jour Ann. of Math.
\vol 120
\yr 1984
\pages 371-397
\endref

\ref\key DJS
\by G. David, J.-L. Journ\'e, and S. Semmes
\paper Operateurs de Calder\'on-Zygmund, fonctions para-accretive
       et interpolation
\jour Rev. Mat. Iberoamericana
\vol 1
\yr 1985
\pages no. 4 1-56
\endref

\ref\key FJ \by M. Frazier and B. Jawerth \paper A discrete
transform and decompositions of distribution spaces \jour J. Funct.
Anal. \vol 93 \yr 1990 \pages 34-170
\endref

\ref\key FS \by C. Fefferman and E. M. Stein \paper Some maximal
inequalities \jour Amer. J. Math. \vol 93 \yr 1971 \pages 107-115
\endref

\ref\key H
\by Y. Han
\paper Calder\'on-type reproducing formula and the $Tb$ theorem
\jour Rev. Mat. Iberoamericana
\vol 10
\yr 1994
\pages 51-91
\endref


\ref\key HLLL
\by Y. Han, M.-Y. Lee, C.-C. Lin, and Y.-C. Lin
\paper Calder\'on-Zygmund operators on product Hardy spaces
\jour J. Funct. Anal.
\vol 258
\yr 2010
\pages 2834-2861
\endref

\ref\key HM
\by T. Hyt\"onen and H. Martikainen
\paper Non-homogeneous $T1$ theorem for bi-parameter singular integrals
\jour preprint available at http://arxiv.org/abs/1209.4473
\endref

\ref\key LZ
\by R.-L. Long and X.-X. Zhu
\paper $L^2$-boundedness of some singular integral operators on product domains
\jour Sci. China Ser. A
\vol 36
\yr 1993
\pages 538-549
\endref

\ref\key J
\by J.-L. Journ\'e
\paper Calder\'on-Zygmund operators on product spaces
\jour Rev. Mat. Iberoamericana
\vol 1
\yr 1985
\pages 55-91
\endref

\ref\key Ma
\by H. Martikainen
\paper Representation of bi-parameter singular integrals by dyadic operators
\jour preprint available at http://arxiv.org/abs/1110.1890
\endref

\ref\key MM
\by A. McIntosh and Y. Meyer
\paper Alg\`ebres d'op\'erateurs d\'efinis par int\'egrales singuli\`eres
\jour C. R. Acad. Sci. Paris S\'er. I Math.
\vol 301
\yr 1985
\pages 395-397
\endref

\ref\key O
\by Y. Ou
\paper A $T(b)$ theorem on product spaces
\jour preprint available at http://arxiv.org/abs/1305.1691
\endref

\ref\key PV
\by S. Pott and P. Villarroya
\paper A $T(1)$ theorem on product spaces
\jour preprint available at http://arxiv.org/abs/\break1105.2516
\endref


\endRefs
\enddocument